\documentclass{article}

\usepackage{graphicx}

\usepackage[french,english]{babel}

\usepackage[utf8]{inputenc}
\usepackage{amsmath, amssymb, amsthm}
\usepackage[alphabetic]{amsrefs}

\usepackage{subfiles}
\usepackage{xr} 
\usepackage{hyperref}
\usepackage[capitalise]{cleveref}
\usepackage{enumitem}

\usepackage{mathtools}

\usepackage{mathrsfs}
\usepackage[scr=rsfso]{mathalpha}

\usepackage{geometry}
\newcommand{\op}[1]{\operatorname{#1}}

\newcommand{\N}{\mathbb{N}}
\newcommand{\R}{\mathbb{R}}
\newcommand{\Z}{\mathbb{Z}}

\newcommand{\<}{\langle}
\renewcommand{\>}{\rangle}

\newcommand{\W}{\mathcal{W}}
\newcommand{\V}{\mathcal{V}}

\usepackage{subcaption}
\usepackage{graphicx}

\newtheorem{theorem}{Theorem}[section]
\newtheorem{proposition}[theorem]{Proposition}

\newtheorem{corollary}[theorem]{Corollary}
\newtheorem{lemma}[theorem]{Lemma}

\theoremstyle{definition}
\newtheorem{definition}[theorem]{Definition}
\newtheorem{example}[theorem]{Example}
\newtheorem{remark}[theorem]{Remark}

\usepackage{tikz}
\usetikzlibrary{shapes.geometric}
\usetikzlibrary{calc}
\usetikzlibrary{quotes, angles}
\usetikzlibrary{graphs}
\usetikzlibrary{positioning}
\usetikzlibrary{3d}
\usetikzlibrary{intersections}
\usetikzlibrary{decorations.markings, decorations.pathmorphing}

\usepackage{newpxtext}
	
\usepackage{authblk}
\title{CAT(0) Coxeter cell complexes, Bipodality and the Parallel Wall Property}
\author{Carl Kristof-Tessier\footnote{Supported by the FRQNT Doctoral Research Scholarship}}
\date{}
\affil{Department of Mathematics and Statistics, McGill University, Burnside Hall, 805 Sherbrooke Street West, Montreal, QC, H3A 0B9, Canada}

\begin{document}
\maketitle

\begin{abstract}
    Coxeter cell complexes are a simultaneous generalization to cube complexes and the Davis complex of a Coxeter group. We show that 1-skeleta of CAT(0) Coxeter cell complexes isometrically embed into a CAT(0) cube complex. When the resulting cube complex is finite dimensional, we show that the condition of bipodality (satisfied by both CAT(0) cube complexes and Davis complexes of Coxeter groups) implies local finiteness of the cubulation. To this end, we show a Parallel Wall Theorem for bipodal CAT(0) Coxeter cell complexes satisfying a uniform bound on the size of a collection of pairwise intersecting walls. In particular, our methods give a new, purely geometric proof of the Parallel Wall Theorem for Coxeter groups. 
\end{abstract}

\section{Introduction}

A \textbf{Coxeter cell complex} is a cell complex built out of polytopes whose symmetry groups are finite Coxeter groups. More precisely, Coxeter cells are Davis complexes of finite Coxeter groups, in the sense of \cite{Davis2008GeomCoxeterGps}. We say a Coxeter cell complex is \textbf{CAT(0)} if it is built out of finitely many types of Coxeter cells and its induced piecewise Euclidean metric is a CAT(0) metric. Motivating examples of such complexes are finite dimensional cube complexes, whose $n$-cells correspond the Davis complex of the Coxeter group $(\Z/2\Z)^n$, and the Davis complex of any Coxeter group, whose cells are formed by the Coxeter cells associated to the group's finite parabolic subgroups. 
As in these examples, CAT(0) Coxeter cell complexes come equipped with a natural wallspace structure. For Coxeter groups, this wallspace structure has been thoroughly studied to deduce combinatorial properties of the groups. At the heart of this study lie two key properties:

\vspace{0.5em}
\textbf{The BPI Property} (Bounded Pairwise Intersections): There is a bound, depending only on the complex, on the number of pairwise intersecting walls.
\vspace{0.5em}

\textbf{The PW Property} (Parallel Wall Property):  There is a constant $C$, depending only on the complex, such that for any wall $\W$ and vertex $x$ at distance $\geq C$ from $\W$, there is a wall $\W'$ separating~$x$ from $\W$.
\vspace{0.5em}

These two properties were first shown to hold for Coxeter groups by Brink and Howlett in \cite{BrinkHowlett1993}.
In general, the PW Property is difficult to establish. We have shown the PW Property always holds in the 2-dimensional case in \cite{Kristof-Tessier25}. The question is currently open in higher dimensions. Here, we give a criterion allowing us to guarantee the PW Property.

We say that $X$ is \textbf{bipodal} if it satisfies the following: Suppose that $\W, \V$ are intersecting walls sharing a common adjacent vertex $y$, and that $x$ is a vertex not separated from $y$ by $\W$ and $\V$. If~$\W$ is separated from $x$ by a wall $\W'$, then for any wall $\mathcal{U}\neq \V$ containing $\W\cap\V$, there is a wall $\mathcal{U}'$ separating $x$ from $\mathcal{U}$.
It is shown by Dyer and Hohwleg in \cite{dyer2016smallrootslowelements} that the Davis complex of a Coxeter group is bipodal (see also \cite{osajdaPrzytycki2022coxeter} for an elementary proof). 

\begin{theorem}\label{thm:bip_BPI_implies_PW}
    Suppose that $X$ is a bipodal $\mathrm{CAT(0)}$ Coxeter cell complex satisfying the BPI Property. Then $X$ satisfies the PW Property. 
\end{theorem}

We emphasize that the original proof of the Parallel Wall theorem for Coxeter groups makes extensive use of the algebraic language coming from root systems. As such, the proof of \cref{thm:bip_BPI_implies_PW} gives a new, geometric proof of the Parallel Wall theorem. 

Using the BPI and PW Properties, it can be deduced from \cite{NibloReeves2003} and \cite{Caprace06} that Coxeter groups act properly on finite dimensional, locally finite CAT(0) cube complexes. In general, Sageev's construction allows us to isometrically embed the 1-skeleton of a CAT(0) Coxeter cell complex $X$ into a cube complex $C(X)$, whose $n$-cubes correspond to collections of $n$ pairwise intersecting walls in $X$. 
In \cite{HruskaWise14}*{Thm 5.7}, Hruska and Wise show that if the complex is itself locally finite, then local finiteness of $C(X)$ is equivalent to $X$ having \textit{compact-wall separation}, an analog of the PW Property for compact sets. We show that for locally finite CAT(0) Coxeter cell complexes, the PW Property implies the compact-wall separation (\cref{cor:loc_fn_cubulations_2}). In particular, we obtain the following:

\begin{corollary}\label[corollary]{cor:loc_fn_cubulations}
    Suppose that $X$ is a bipodal locally finite $\mathrm{CAT(0)}$ Coxeter cell complex satisfying the BPI Property. Then $X^{(1)}$ isometrically embeds into a finite dimensional, locally finite $\mathrm{CAT(0)}$ cube complex.
\end{corollary}

Furthermore, metrically proper group actions on $X$ extend to metrically proper actions on~$C(X)$ \cite{HruskaWise14}. It follows, for example, that if $G$ is a group acting properly on $X$, and $G$ has a bound on the size of its finite subgroups (e.g. when the action is cocompact or, more generally, when $X$ is uniformly locally finite), then $G$ satisfies the Tits alternative \cite{SageevWise}.

\begin{corollary}\label[corollary]{cor:proper_action}
    If $G$ acts properly by automorphisms on a bipodal locally finite $\mathrm{CAT(0)}$ Coxeter cell complex satisfying the BPI Property, then $G$ acts properly on a locally finite, finite dimensional $\mathrm{CAT(0)}$ cube complex.
\end{corollary}

Contrary to much of the theory of CAT(0) cube complexes, we rely on the CAT(0) metric in most of our proofs. Basic results in CAT(0) geometry, as well as background on Coxeter cell complexes, can be found in \cref{sec:background}.
In \cref{sec:Sageev_construction}, we recall Sageev's construction and show Corollary \ref{cor:loc_fn_cubulations} under the assumption of \cref{thm:bip_BPI_implies_PW}.
In \cref{sec:angles}, we introduce the notion of \textit{angle between walls}, and use this to give a CAT(0)-like inequality for ``wall triangles''. We also generalize \cite{BrinkHowlett1993}*{Lem~1.7}, showing that these angles can be used to compute the distance between a vertex and a wall. 
Finally, we prove \cref{thm:bip_BPI_implies_PW} in \cref{sec:proof}. 

\vspace{1em}
\textbf{Acknowledgements}. The author would like to thank Piotr Przytycki for his helpful discussions throughout this project, and for reviewing preliminary versions of this paper. The author would also like to thank Dani Wise for providing initial feedback on the author's thesis, which served as the foundation for this work. We acknowledge that AI was not used in any capacity in the process of writing this paper, or proving the results therein.

\section{Background}\label{sec:background}

\subsection{CAT(0) Geometry}
Let $(X,d)$ be a geodesic metric space. We denote by $[x,y]$ a geodesic segment connecting $x$ and~$y$. A \textbf{geodesic triangle} in $X$ consists of a triple of points $x,y,z\in X$, together with geodesic segments between the three points. We denote by $\Delta(x,y,z)$ the union of these segments. 
For every geodesic triangle in $X$ with vertices $x,y,z\in X$, there is a unique (up to isometry) Euclidean triangle $\overline{\Delta}(x,y,z)$ in $\R^2$, called the \textbf{comparison triangle}, with vertices $\overline x, \overline y,\overline z\in \R^2$ and side lengths $d_{\R^2}(\overline x,\overline y)=d(x,y),d_{\R^2}(\overline x,\overline z)=d(x,z)$ and $d_{\R^2}(\overline y,\overline z)=d(y,z)$. A point $\overline p\in [\overline x, \overline y]$ (resp $[\overline x, \overline z]$ and $[\overline y, \overline z]$) is called a \textbf{comparison point} for $p\in [x,y]$ (resp $[x, z]$ and $[y, z]$) if $d_{\R^2}(\overline x, \overline p) = d(x, p)$.

\begin{definition}
    A geodesic space $X$ is \textbf{CAT(0)} if for any geodesic triangle $\Delta$ and comparison points $\overline p, \overline q\in \overline{\Delta}$ for $p,q\in \Delta$, we have $d(p,q)\leq d_{\R^2}(\bar{p},\bar{q})$. We say that $X$ is \textbf{non-positively curved} if it is locally CAT(0).
\end{definition}

It is a standard exercise to show that CAT(0) spaces are uniquely geodesic and thus contractible. Conversely, the Cartan-Hadamard theorem states that if $X$ is complete, non-positively curved and simply connected, then $X$ is CAT(0) \cite{Bridson1999}*{Thm II.4.1}.

There is also an equivalent definition for CAT(0) spaces in terms of angles. Given $x,y,z\in X$, the \textbf{comparison angle} $\overline{\angle_z}(x,y)$ is the Euclidean angle in the comparison triangle $\overline{\Delta}(x,y,z)$ at the vertex $\overline{z}$.
Given geodesics $\gamma_1, \gamma_2$ with $\gamma_1(0)=\gamma_2(0)=p$, we may then define the \textbf{Alexandrov angle} between $\gamma_1$ and $\gamma_2$ as
\[
        \angle(\gamma_1,\gamma_2)=\lim_{\epsilon\to 0} \sup_{t,t'\leq \epsilon} \overline{\angle_p}(\gamma_1(t),\gamma_2(t'))
\]
We denote by $\angle_z(x,y)$ the Alexandrov angle between segments $[z,x]$ and $[z,y]$.

\begin{proposition}[{\cite[Prop II.1.7]{Bridson1999}}]
    \label[proposition]{prop:cat0_angle_ineq}
    $X$ is $\mathrm{CAT(0)}$ if and only if $\angle_z(x,y)\leq \overline{\angle_z}(x,y)$ for all $x,y,z\in X$. In particular, the sum of Alexandrov angles in a geodesic triangle of a $\mathrm{CAT(0)}$ space is bounded by $\pi$.
\end{proposition}

The metric of a CAT(0) space is convex, and due to this CAT(0) spaces behave nicely with respect to their convex subsets. 

\begin{proposition}[{\cite[Prop II.2.4]{Bridson1999}}]
\label[proposition]{prop:projection_existence_uniqueness}
    Let $X$ be a complete $\mathrm{CAT(0)}$ space, and $C\subset X$ a closed convex subset of $X$. Then for $x\in X$, there is a unique point $\pi_C(x)\in C$ such that $d(x,\pi_C(x))=d(x,C)$. Furthermore, if $x\notin C$ then $\pi_C(x)$ is the unique point in $X$ satisfying $\angle_{\pi_C(x)}(x,y)\geq \pi/2$ for all $y\in C\setminus\{\pi_C(x)\}$.
\end{proposition}

As a consequence, one can show that the function $x\to d(x,C)$ is a convex function \cite[Cor~II.2.5]{Bridson1999}, and therefore the $r$-\textbf{neighborhood} of $C$, defined by $N_r(C)= \{x\in X:d(x,C)< r\}$, is convex in $X$ for any $r\geq 0$.

\subsection{Coxeter groups and Coxeter cell complexes}

Let $\Gamma$ be a simplicial graph with edge labels taking values in $\{2,3,\dots\}$. Then the \textbf{Coxeter group} with defining graph $\Gamma$ is
    \[
        W_\Gamma = \left\< V(\Gamma) \Bigm| s^2 =1, (st)^{m_{st}}=1 \text{ for each}
        \begin{matrix}\begin{tikzpicture}
            \filldraw (0, 0) circle (1pt) node [left] {$s$};
            \filldraw (.75, 0) circle (1pt) node [right] {$t$};
            \draw (0, 0) -- (.75, 0) node[midway, above] {\tiny $m_{st}$};
        \end{tikzpicture}\end{matrix}
        \in E(\Gamma) \right\>
    \]

A Coxeter group $W_\Gamma$ has an associated faithful representation called the \textbf{Tits representation}. We consider a vector space $V$ with basis $\{\alpha_s\}_{s\in V(\Gamma)}$, equipped with the symmetric bilinear form $\<\alpha_s, \alpha_t\> = -\cos\frac{\pi}{m_{st}}$.
Then $W_\Gamma$ acts faithfully on $V$ via 
\[
    s\cdot v = v - 2\<\alpha_s, v\> \alpha_s
\]
for $s\in V(\Gamma)$ \cite[Sec 6.12]{Davis2008GeomCoxeterGps}. 
A Coxeter group $W_\Gamma$ is finite if and only if the bilinear form arising from the Tits representation is positive definite \cite[Thm 6.12.9]{Davis2008GeomCoxeterGps}. In this case, the representation gives an action of $W_\Gamma$ by isometries on a Euclidean space $V$, where conjugates of elements in $V(\Gamma)$ act by reflections along hyperplanes in $V$. These hyperplanes partition $V$ into a set of cones, on which $W_\Gamma$ acts transitively. The \textbf{fundamental cone in $V$} is the cone supported by the hyperplanes fixed by the generators of $W_\Gamma$. We denote by $e$ the point at distance $\frac{1}{2}$ from each supporting hyperplane in the fundamental cone.

\begin{definition}
    The \textbf{Coxeter cell} of type $\Gamma$, denoted by $Z_\Gamma$, is the convex hull of $W_\Gamma\cdot e$. 
\end{definition}

The cell structure on $Z_\Gamma$ can be described in terms of the Coxeter subgroups $W_{\Gamma'}\leq W_\Gamma$ generated by subgraphs $\Gamma'\subset \Gamma$. 
As shown in  \cite[Lem 7.3.3]{Davis2008GeomCoxeterGps}, every codimension-1 face of a Coxeter cell $Z_\Gamma$ corresponds to a coset  $gW_{\Gamma'}$, where $\Gamma'\subset \Gamma$ is a full subgraph with size $\#V(\Gamma)-1$. The face corresponding to $gW_{\Gamma'}$ is isometric to the convex hull of $gW_{\Gamma'}\cdot e$. Thus each $Z_{\Gamma'}$ embeds into~$Z_\Gamma$, allowing us to inductively define the cell structure on $Z_\Gamma$.  In particular, the vertices of $Z_\Gamma$ correspond to elements of $W_\Gamma$, and we have a clear vertex-transitive action of $W_\Gamma$ on $Z_\Gamma$. 

\begin{figure}
    \centering
    \begin{subfigure}{0.25\textwidth}
        \centering
        \begin{tikzpicture}
                \foreach \k in {0,1,...,5}
                \draw ({30+(\k-1)*60}:1) -- ({30+\k*60}:1);
        \end{tikzpicture}
        
        \begin{minipage}{2cm}
            \vfill
        \end{minipage}

        \begin{minipage}{0.1cm}
            \vfill
        \end{minipage}
    \end{subfigure}
    \begin{subfigure}{0.25\textwidth}
        \centering
        \begin{tikzpicture}
            	[x={(-0.702577cm, -0.347995cm)},
            	y={(0.711608cm, -0.343626cm)},
            	z={(0.000039cm, 0.872250cm)},
            	scale=1.000000,
            	back/.style={loosely dotted, thin, opacity=0},
            	edge/.style={color=black!95!black, thick},
            	facet/.style={fill=gray,fill opacity=0.000000},
            	vertex/.style={}]
            %
            %
            \coordinate (1.00000, -1.00000, -1.00000) at (1.00000, -1.00000, -1.00000);
            \coordinate (1.00000, 1.00000, -1.00000) at (1.00000, 1.00000, -1.00000);
            \coordinate (1.00000, 1.00000, 1.00000) at (1.00000, 1.00000, 1.00000);
            \coordinate (1.00000, -1.00000, 1.00000) at (1.00000, -1.00000, 1.00000);
            \coordinate (-1.00000, -1.00000, 1.00000) at (-1.00000, -1.00000, 1.00000);
            \coordinate (-1.00000, -1.00000, -1.00000) at (-1.00000, -1.00000, -1.00000);
            \coordinate (-1.00000, 1.00000, -1.00000) at (-1.00000, 1.00000, -1.00000);
            \coordinate (-1.00000, 1.00000, 1.00000) at (-1.00000, 1.00000, 1.00000);
            \draw[edge,back] (1.00000, -1.00000, -1.00000) -- (-1.00000, -1.00000, -1.00000);
            \draw[edge,back] (-1.00000, -1.00000, 1.00000) -- (-1.00000, -1.00000, -1.00000);
            \draw[edge,back] (-1.00000, -1.00000, -1.00000) -- (-1.00000, 1.00000, -1.00000);
            \node[vertex] at (-1.00000, -1.00000, -1.00000)     {};
            \fill[facet] (1.00000, -1.00000, 1.00000) -- (1.00000, -1.00000, -1.00000) -- (1.00000, 1.00000, -1.00000) -- (1.00000, 1.00000, 1.00000) -- cycle {};
            \fill[facet] (-1.00000, 1.00000, 1.00000) -- (1.00000, 1.00000, 1.00000) -- (1.00000, 1.00000, -1.00000) -- (-1.00000, 1.00000, -1.00000) -- cycle {};
            \fill[facet] (-1.00000, 1.00000, 1.00000) -- (1.00000, 1.00000, 1.00000) -- (1.00000, -1.00000, 1.00000) -- (-1.00000, -1.00000, 1.00000) -- cycle {};
            \draw[edge] (1.00000, -1.00000, -1.00000) -- (1.00000, 1.00000, -1.00000);
            \draw[edge] (1.00000, -1.00000, -1.00000) -- (1.00000, -1.00000, 1.00000);
            \draw[edge] (1.00000, 1.00000, -1.00000) -- (1.00000, 1.00000, 1.00000);
            \draw[edge] (1.00000, 1.00000, -1.00000) -- (-1.00000, 1.00000, -1.00000);
            \draw[edge] (1.00000, 1.00000, 1.00000) -- (1.00000, -1.00000, 1.00000);
            \draw[edge] (1.00000, 1.00000, 1.00000) -- (-1.00000, 1.00000, 1.00000);
            \draw[edge] (1.00000, -1.00000, 1.00000) -- (-1.00000, -1.00000, 1.00000);
            \draw[edge] (-1.00000, -1.00000, 1.00000) -- (-1.00000, 1.00000, 1.00000);
            \draw[edge] (-1.00000, 1.00000, -1.00000) -- (-1.00000, 1.00000, 1.00000);
            \node[vertex] at (1.00000, -1.00000, -1.00000)     {};
            \node[vertex] at (1.00000, 1.00000, -1.00000)     {};
            \node[vertex] at (1.00000, 1.00000, 1.00000)     {};
            \node[vertex] at (1.00000, -1.00000, 1.00000)     {};
            \node[vertex] at (-1.00000, -1.00000, 1.00000)     {};
            \node[vertex] at (-1.00000, 1.00000, -1.00000)     {};
            \node[vertex] at (-1.00000, 1.00000, 1.00000)     {};
            \end{tikzpicture}
            
            \begin{minipage}{0.1cm}
            \vfill
        \end{minipage}
    \end{subfigure}
    \begin{subfigure}{0.40\textwidth}
        \centering
        \begin{tikzpicture}%
            	[x={(0.249656cm, -0.577639cm)},
            	y={(0.777700cm, -0.358578cm)},
            	z={(-0.576936cm, -0.733318cm)},
            	scale=0.500000,
            	back/.style={loosely dotted, thin, color=white},
            	edge/.style={color=black!95!black, thick},
            	facet/.style={fill=gray,fill opacity=0.000000},
            	vertex/.style={}]
            %
            %
            \coordinate (2.41421, 1.00000, 3.82843) at (2.41421, 1.00000, 3.82843);
            \coordinate (1.00000, 2.41421, 3.82843) at (1.00000, 2.41421, 3.82843);
            \coordinate (1.00000, 3.82843, 2.41421) at (1.00000, 3.82843, 2.41421);
            \coordinate (2.41421, 1.00000, -3.82843) at (2.41421, 1.00000, -3.82843);
            \coordinate (1.00000, 2.41421, -3.82843) at (1.00000, 2.41421, -3.82843);
            \coordinate (1.00000, 3.82843, -2.41421) at (1.00000, 3.82843, -2.41421);
            \coordinate (2.41421, 3.82843, 1.00000) at (2.41421, 3.82843, 1.00000);
            \coordinate (2.41421, 3.82843, -1.00000) at (2.41421, 3.82843, -1.00000);
            \coordinate (1.00000, -2.41421, 3.82843) at (1.00000, -2.41421, 3.82843);
            \coordinate (1.00000, -3.82843, 2.41421) at (1.00000, -3.82843, 2.41421);
            \coordinate (2.41421, -1.00000, 3.82843) at (2.41421, -1.00000, 3.82843);
            \coordinate (1.00000, -2.41421, -3.82843) at (1.00000, -2.41421, -3.82843);
            \coordinate (1.00000, -3.82843, -2.41421) at (1.00000, -3.82843, -2.41421);
            \coordinate (2.41421, -1.00000, -3.82843) at (2.41421, -1.00000, -3.82843);
            \coordinate (3.82843, 1.00000, 2.41421) at (3.82843, 1.00000, 2.41421);
            \coordinate (3.82843, 1.00000, -2.41421) at (3.82843, 1.00000, -2.41421);
            \coordinate (3.82843, 2.41421, 1.00000) at (3.82843, 2.41421, 1.00000);
            \coordinate (3.82843, 2.41421, -1.00000) at (3.82843, 2.41421, -1.00000);
            \coordinate (3.82843, -1.00000, 2.41421) at (3.82843, -1.00000, 2.41421);
            \coordinate (3.82843, -1.00000, -2.41421) at (3.82843, -1.00000, -2.41421);
            \coordinate (2.41421, -3.82843, 1.00000) at (2.41421, -3.82843, 1.00000);
            \coordinate (2.41421, -3.82843, -1.00000) at (2.41421, -3.82843, -1.00000);
            \coordinate (3.82843, -2.41421, 1.00000) at (3.82843, -2.41421, 1.00000);
            \coordinate (3.82843, -2.41421, -1.00000) at (3.82843, -2.41421, -1.00000);
            \coordinate (-2.41421, 1.00000, 3.82843) at (-2.41421, 1.00000, 3.82843);
            \coordinate (-1.00000, 2.41421, 3.82843) at (-1.00000, 2.41421, 3.82843);
            \coordinate (-1.00000, 3.82843, 2.41421) at (-1.00000, 3.82843, 2.41421);
            \coordinate (-2.41421, 1.00000, -3.82843) at (-2.41421, 1.00000, -3.82843);
            \coordinate (-1.00000, 2.41421, -3.82843) at (-1.00000, 2.41421, -3.82843);
            \coordinate (-1.00000, 3.82843, -2.41421) at (-1.00000, 3.82843, -2.41421);
            \coordinate (-2.41421, 3.82843, 1.00000) at (-2.41421, 3.82843, 1.00000);
            \coordinate (-2.41421, 3.82843, -1.00000) at (-2.41421, 3.82843, -1.00000);
            \coordinate (-1.00000, -2.41421, 3.82843) at (-1.00000, -2.41421, 3.82843);
            \coordinate (-1.00000, -3.82843, 2.41421) at (-1.00000, -3.82843, 2.41421);
            \coordinate (-2.41421, -1.00000, 3.82843) at (-2.41421, -1.00000, 3.82843);
            \coordinate (-1.00000, -2.41421, -3.82843) at (-1.00000, -2.41421, -3.82843);
            \coordinate (-1.00000, -3.82843, -2.41421) at (-1.00000, -3.82843, -2.41421);
            \coordinate (-2.41421, -1.00000, -3.82843) at (-2.41421, -1.00000, -3.82843);
            \coordinate (-3.82843, 1.00000, 2.41421) at (-3.82843, 1.00000, 2.41421);
            \coordinate (-3.82843, 1.00000, -2.41421) at (-3.82843, 1.00000, -2.41421);
            \coordinate (-3.82843, 2.41421, 1.00000) at (-3.82843, 2.41421, 1.00000);
            \coordinate (-3.82843, 2.41421, -1.00000) at (-3.82843, 2.41421, -1.00000);
            \coordinate (-3.82843, -1.00000, 2.41421) at (-3.82843, -1.00000, 2.41421);
            \coordinate (-3.82843, -1.00000, -2.41421) at (-3.82843, -1.00000, -2.41421);
            \coordinate (-2.41421, -3.82843, 1.00000) at (-2.41421, -3.82843, 1.00000);
            \coordinate (-2.41421, -3.82843, -1.00000) at (-2.41421, -3.82843, -1.00000);
            \coordinate (-3.82843, -2.41421, 1.00000) at (-3.82843, -2.41421, 1.00000);
            \coordinate (-3.82843, -2.41421, -1.00000) at (-3.82843, -2.41421, -1.00000);
            \draw[edge,back] (2.41421, 1.00000, -3.82843) -- (1.00000, 2.41421, -3.82843);
            \draw[edge,back] (2.41421, 1.00000, -3.82843) -- (2.41421, -1.00000, -3.82843);
            \draw[edge,back] (2.41421, 1.00000, -3.82843) -- (3.82843, 1.00000, -2.41421);
            \draw[edge,back] (2.41421, 3.82843, -1.00000) -- (3.82843, 2.41421, -1.00000);
            \draw[edge,back] (1.00000, -2.41421, 3.82843) -- (1.00000, -3.82843, 2.41421);
            \draw[edge,back] (1.00000, -3.82843, 2.41421) -- (2.41421, -3.82843, 1.00000);
            \draw[edge,back] (1.00000, -3.82843, 2.41421) -- (-1.00000, -3.82843, 2.41421);
            \draw[edge,back] (2.41421, -1.00000, 3.82843) -- (3.82843, -1.00000, 2.41421);
            \draw[edge,back] (1.00000, -2.41421, -3.82843) -- (1.00000, -3.82843, -2.41421);
            \draw[edge,back] (1.00000, -2.41421, -3.82843) -- (2.41421, -1.00000, -3.82843);
            \draw[edge,back] (1.00000, -2.41421, -3.82843) -- (-1.00000, -2.41421, -3.82843);
            \draw[edge,back] (1.00000, -3.82843, -2.41421) -- (2.41421, -3.82843, -1.00000);
            \draw[edge,back] (1.00000, -3.82843, -2.41421) -- (-1.00000, -3.82843, -2.41421);
            \draw[edge,back] (2.41421, -1.00000, -3.82843) -- (3.82843, -1.00000, -2.41421);
            \draw[edge,back] (3.82843, 1.00000, 2.41421) -- (3.82843, -1.00000, 2.41421);
            \draw[edge,back] (3.82843, 1.00000, -2.41421) -- (3.82843, 2.41421, -1.00000);
            \draw[edge,back] (3.82843, 1.00000, -2.41421) -- (3.82843, -1.00000, -2.41421);
            \draw[edge,back] (3.82843, 2.41421, 1.00000) -- (3.82843, 2.41421, -1.00000);
            \draw[edge,back] (3.82843, -1.00000, 2.41421) -- (3.82843, -2.41421, 1.00000);
            \draw[edge,back] (3.82843, -1.00000, -2.41421) -- (3.82843, -2.41421, -1.00000);
            \draw[edge,back] (2.41421, -3.82843, 1.00000) -- (2.41421, -3.82843, -1.00000);
            \draw[edge,back] (2.41421, -3.82843, 1.00000) -- (3.82843, -2.41421, 1.00000);
            \draw[edge,back] (2.41421, -3.82843, -1.00000) -- (3.82843, -2.41421, -1.00000);
            \draw[edge,back] (3.82843, -2.41421, 1.00000) -- (3.82843, -2.41421, -1.00000);
            \draw[edge,back] (-1.00000, -2.41421, -3.82843) -- (-1.00000, -3.82843, -2.41421);
            \draw[edge,back] (-1.00000, -2.41421, -3.82843) -- (-2.41421, -1.00000, -3.82843);
            \draw[edge,back] (-1.00000, -3.82843, -2.41421) -- (-2.41421, -3.82843, -1.00000);
            \node[vertex] at (1.00000, -3.82843, 2.41421)     {};
            \node[vertex] at (1.00000, -2.41421, -3.82843)     {};
            \node[vertex] at (1.00000, -3.82843, -2.41421)     {};
            \node[vertex] at (-1.00000, -2.41421, -3.82843)     {};
            \node[vertex] at (-1.00000, -3.82843, -2.41421)     {};
            \node[vertex] at (2.41421, 1.00000, -3.82843)     {};
            \node[vertex] at (2.41421, -1.00000, -3.82843)     {};
            \node[vertex] at (2.41421, -3.82843, 1.00000)     {};
            \node[vertex] at (2.41421, -3.82843, -1.00000)     {};
            \node[vertex] at (3.82843, -1.00000, 2.41421)     {};
            \node[vertex] at (3.82843, 1.00000, -2.41421)     {};
            \node[vertex] at (3.82843, -1.00000, -2.41421)     {};
            \node[vertex] at (3.82843, 2.41421, -1.00000)     {};
            \node[vertex] at (3.82843, -2.41421, 1.00000)     {};
            \node[vertex] at (3.82843, -2.41421, -1.00000)     {};
            \fill[facet] (-1.00000, 3.82843, 2.41421) -- (1.00000, 3.82843, 2.41421) -- (1.00000, 2.41421, 3.82843) -- (-1.00000, 2.41421, 3.82843) -- cycle {};
            \fill[facet] (-1.00000, 3.82843, -2.41421) -- (1.00000, 3.82843, -2.41421) -- (1.00000, 2.41421, -3.82843) -- (-1.00000, 2.41421, -3.82843) -- cycle {};
            \fill[facet] (-2.41421, -1.00000, 3.82843) -- (-2.41421, 1.00000, 3.82843) -- (-1.00000, 2.41421, 3.82843) -- (1.00000, 2.41421, 3.82843) -- (2.41421, 1.00000, 3.82843) -- (2.41421, -1.00000, 3.82843) -- (1.00000, -2.41421, 3.82843) -- (-1.00000, -2.41421, 3.82843) -- cycle {};
            \fill[facet] (-2.41421, 3.82843, -1.00000) -- (-1.00000, 3.82843, -2.41421) -- (1.00000, 3.82843, -2.41421) -- (2.41421, 3.82843, -1.00000) -- (2.41421, 3.82843, 1.00000) -- (1.00000, 3.82843, 2.41421) -- (-1.00000, 3.82843, 2.41421) -- (-2.41421, 3.82843, 1.00000) -- cycle {};
            \fill[facet] (-3.82843, -1.00000, 2.41421) -- (-2.41421, -1.00000, 3.82843) -- (-2.41421, 1.00000, 3.82843) -- (-3.82843, 1.00000, 2.41421) -- cycle {};
            \fill[facet] (-3.82843, -1.00000, -2.41421) -- (-2.41421, -1.00000, -3.82843) -- (-2.41421, 1.00000, -3.82843) -- (-3.82843, 1.00000, -2.41421) -- cycle {};
            \fill[facet] (3.82843, 2.41421, 1.00000) -- (2.41421, 3.82843, 1.00000) -- (1.00000, 3.82843, 2.41421) -- (1.00000, 2.41421, 3.82843) -- (2.41421, 1.00000, 3.82843) -- (3.82843, 1.00000, 2.41421) -- cycle {};
            \fill[facet] (-3.82843, 2.41421, 1.00000) -- (-2.41421, 3.82843, 1.00000) -- (-1.00000, 3.82843, 2.41421) -- (-1.00000, 2.41421, 3.82843) -- (-2.41421, 1.00000, 3.82843) -- (-3.82843, 1.00000, 2.41421) -- cycle {};
            \fill[facet] (-3.82843, 2.41421, -1.00000) -- (-2.41421, 3.82843, -1.00000) -- (-1.00000, 3.82843, -2.41421) -- (-1.00000, 2.41421, -3.82843) -- (-2.41421, 1.00000, -3.82843) -- (-3.82843, 1.00000, -2.41421) -- cycle {};
            \fill[facet] (-3.82843, 2.41421, -1.00000) -- (-2.41421, 3.82843, -1.00000) -- (-2.41421, 3.82843, 1.00000) -- (-3.82843, 2.41421, 1.00000) -- cycle {};
            \fill[facet] (-3.82843, -2.41421, 1.00000) -- (-3.82843, -1.00000, 2.41421) -- (-2.41421, -1.00000, 3.82843) -- (-1.00000, -2.41421, 3.82843) -- (-1.00000, -3.82843, 2.41421) -- (-2.41421, -3.82843, 1.00000) -- cycle {};
            \fill[facet] (-3.82843, -2.41421, -1.00000) -- (-3.82843, -1.00000, -2.41421) -- (-3.82843, 1.00000, -2.41421) -- (-3.82843, 2.41421, -1.00000) -- (-3.82843, 2.41421, 1.00000) -- (-3.82843, 1.00000, 2.41421) -- (-3.82843, -1.00000, 2.41421) -- (-3.82843, -2.41421, 1.00000) -- cycle {};
            \fill[facet] (-3.82843, -2.41421, -1.00000) -- (-2.41421, -3.82843, -1.00000) -- (-2.41421, -3.82843, 1.00000) -- (-3.82843, -2.41421, 1.00000) -- cycle {};
            \draw[edge] (2.41421, 1.00000, 3.82843) -- (1.00000, 2.41421, 3.82843);
            \draw[edge] (2.41421, 1.00000, 3.82843) -- (2.41421, -1.00000, 3.82843);
            \draw[edge] (2.41421, 1.00000, 3.82843) -- (3.82843, 1.00000, 2.41421);
            \draw[edge] (1.00000, 2.41421, 3.82843) -- (1.00000, 3.82843, 2.41421);
            \draw[edge] (1.00000, 2.41421, 3.82843) -- (-1.00000, 2.41421, 3.82843);
            \draw[edge] (1.00000, 3.82843, 2.41421) -- (2.41421, 3.82843, 1.00000);
            \draw[edge] (1.00000, 3.82843, 2.41421) -- (-1.00000, 3.82843, 2.41421);
            \draw[edge] (1.00000, 2.41421, -3.82843) -- (1.00000, 3.82843, -2.41421);
            \draw[edge] (1.00000, 2.41421, -3.82843) -- (-1.00000, 2.41421, -3.82843);
            \draw[edge] (1.00000, 3.82843, -2.41421) -- (2.41421, 3.82843, -1.00000);
            \draw[edge] (1.00000, 3.82843, -2.41421) -- (-1.00000, 3.82843, -2.41421);
            \draw[edge] (2.41421, 3.82843, 1.00000) -- (2.41421, 3.82843, -1.00000);
            \draw[edge] (2.41421, 3.82843, 1.00000) -- (3.82843, 2.41421, 1.00000);
            \draw[edge] (1.00000, -2.41421, 3.82843) -- (2.41421, -1.00000, 3.82843);
            \draw[edge] (1.00000, -2.41421, 3.82843) -- (-1.00000, -2.41421, 3.82843);
            \draw[edge] (3.82843, 1.00000, 2.41421) -- (3.82843, 2.41421, 1.00000);
            \draw[edge] (-2.41421, 1.00000, 3.82843) -- (-1.00000, 2.41421, 3.82843);
            \draw[edge] (-2.41421, 1.00000, 3.82843) -- (-2.41421, -1.00000, 3.82843);
            \draw[edge] (-2.41421, 1.00000, 3.82843) -- (-3.82843, 1.00000, 2.41421);
            \draw[edge] (-1.00000, 2.41421, 3.82843) -- (-1.00000, 3.82843, 2.41421);
            \draw[edge] (-1.00000, 3.82843, 2.41421) -- (-2.41421, 3.82843, 1.00000);
            \draw[edge] (-2.41421, 1.00000, -3.82843) -- (-1.00000, 2.41421, -3.82843);
            \draw[edge] (-2.41421, 1.00000, -3.82843) -- (-2.41421, -1.00000, -3.82843);
            \draw[edge] (-2.41421, 1.00000, -3.82843) -- (-3.82843, 1.00000, -2.41421);
            \draw[edge] (-1.00000, 2.41421, -3.82843) -- (-1.00000, 3.82843, -2.41421);
            \draw[edge] (-1.00000, 3.82843, -2.41421) -- (-2.41421, 3.82843, -1.00000);
            \draw[edge] (-2.41421, 3.82843, 1.00000) -- (-2.41421, 3.82843, -1.00000);
            \draw[edge] (-2.41421, 3.82843, 1.00000) -- (-3.82843, 2.41421, 1.00000);
            \draw[edge] (-2.41421, 3.82843, -1.00000) -- (-3.82843, 2.41421, -1.00000);
            \draw[edge] (-1.00000, -2.41421, 3.82843) -- (-1.00000, -3.82843, 2.41421);
            \draw[edge] (-1.00000, -2.41421, 3.82843) -- (-2.41421, -1.00000, 3.82843);
            \draw[edge] (-1.00000, -3.82843, 2.41421) -- (-2.41421, -3.82843, 1.00000);
            \draw[edge] (-2.41421, -1.00000, 3.82843) -- (-3.82843, -1.00000, 2.41421);
            \draw[edge] (-2.41421, -1.00000, -3.82843) -- (-3.82843, -1.00000, -2.41421);
            \draw[edge] (-3.82843, 1.00000, 2.41421) -- (-3.82843, 2.41421, 1.00000);
            \draw[edge] (-3.82843, 1.00000, 2.41421) -- (-3.82843, -1.00000, 2.41421);
            \draw[edge] (-3.82843, 1.00000, -2.41421) -- (-3.82843, 2.41421, -1.00000);
            \draw[edge] (-3.82843, 1.00000, -2.41421) -- (-3.82843, -1.00000, -2.41421);
            \draw[edge] (-3.82843, 2.41421, 1.00000) -- (-3.82843, 2.41421, -1.00000);
            \draw[edge] (-3.82843, -1.00000, 2.41421) -- (-3.82843, -2.41421, 1.00000);
            \draw[edge] (-3.82843, -1.00000, -2.41421) -- (-3.82843, -2.41421, -1.00000);
            \draw[edge] (-2.41421, -3.82843, 1.00000) -- (-2.41421, -3.82843, -1.00000);
            \draw[edge] (-2.41421, -3.82843, 1.00000) -- (-3.82843, -2.41421, 1.00000);
            \draw[edge] (-2.41421, -3.82843, -1.00000) -- (-3.82843, -2.41421, -1.00000);
            \draw[edge] (-3.82843, -2.41421, 1.00000) -- (-3.82843, -2.41421, -1.00000);
            \node[vertex] at (2.41421, 1.00000, 3.82843)     {};
            \node[vertex] at (1.00000, 2.41421, 3.82843)     {};
            \node[vertex] at (1.00000, 3.82843, 2.41421)     {};
            \node[vertex] at (1.00000, 2.41421, -3.82843)     {};
            \node[vertex] at (1.00000, 3.82843, -2.41421)     {};
            \node[vertex] at (2.41421, 3.82843, 1.00000)     {};
            \node[vertex] at (2.41421, 3.82843, -1.00000)     {};
            \node[vertex] at (1.00000, -2.41421, 3.82843)     {};
            \node[vertex] at (2.41421, -1.00000, 3.82843)     {};
            \node[vertex] at (3.82843, 1.00000, 2.41421)     {};
            \node[vertex] at (3.82843, 2.41421, 1.00000)     {};
            \node[vertex] at (-2.41421, 1.00000, 3.82843)     {};
            \node[vertex] at (-1.00000, 2.41421, 3.82843)     {};
            \node[vertex] at (-1.00000, 3.82843, 2.41421)     {};
            \node[vertex] at (-2.41421, 1.00000, -3.82843)     {};
            \node[vertex] at (-1.00000, 2.41421, -3.82843)     {};
            \node[vertex] at (-1.00000, 3.82843, -2.41421)     {};
            \node[vertex] at (-2.41421, 3.82843, 1.00000)     {};
            \node[vertex] at (-2.41421, 3.82843, -1.00000)     {};
            \node[vertex] at (-1.00000, -2.41421, 3.82843)     {};
            \node[vertex] at (-1.00000, -3.82843, 2.41421)     {};
            \node[vertex] at (-2.41421, -1.00000, 3.82843)     {};
            \node[vertex] at (-2.41421, -1.00000, -3.82843)     {};
            \node[vertex] at (-3.82843, 1.00000, 2.41421)     {};
            \node[vertex] at (-3.82843, 1.00000, -2.41421)     {};
            \node[vertex] at (-3.82843, 2.41421, 1.00000)     {};
            \node[vertex] at (-3.82843, 2.41421, -1.00000)     {};
            \node[vertex] at (-3.82843, -1.00000, 2.41421)     {};
            \node[vertex] at (-3.82843, -1.00000, -2.41421)     {};
            \node[vertex] at (-2.41421, -3.82843, 1.00000)     {};
            \node[vertex] at (-2.41421, -3.82843, -1.00000)     {};
            \node[vertex] at (-3.82843, -2.41421, 1.00000)     {};
            \node[vertex] at (-3.82843, -2.41421, -1.00000)     {};
    \end{tikzpicture}
    \end{subfigure}
    \caption[Examples of Coxeter cells]{Coxeter cells of type \ 
        \protect\begin{tikzpicture}[node distance={20pt}, main/.style = {draw, circle, fill=black, inner sep=1pt}]
           \protect\node[main] (a) {};
           \protect\node[main] (b) [right of=a] {};
           \protect\draw (a) -- (b) node[midway, above] {\small{3}};
        \end{tikzpicture}\protect\ , \         $\begin{matrix}\begin{tikzpicture}[node distance={20pt}, main/.style = {draw, circle, fill=black, inner sep=1pt}]
           \protect\node[main] (a) {};
           \protect\node[main] (b) [above right of=a] {};
           \protect\node[main] (c) [below right of=b] {};
           \protect\draw (a) -- (b) node[midway, above left] {\small{2}};
           \protect\draw (b) -- (c) node[midway, above right] {\small{2}};
           \protect\draw (a) -- (c) node[midway, below] {\small{2}};
        \end{tikzpicture}\end{matrix}$, and \         $\begin{matrix}\begin{tikzpicture}[node distance={20pt}, main/.style = {draw, circle, fill=black, inner sep=1pt}]
           \protect\node[main] (a) {};
           \protect\node[main] (b) [above right of=a] {};
           \protect\node[main] (c) [below right of=b] {};
           \protect\draw (a) -- (b) node[midway, above left] {\small{2}};
           \protect\draw (b) -- (c) node[midway, above right] {\small{3}};
           \protect\draw (a) -- (c) node[midway, below] {\small{4}};
        \end{tikzpicture}\end{matrix} $ (generated using SageMath \cite{sagemath}).}
    \label{fig:coxeter_cells}
\end{figure}

\begin{definition}
    A \textbf{Coxeter cell complex} is a Euclidean polyhedral complex whose cells are isometric to Coxeter cells. A Coxeter cell complex is \textbf{CAT(0)} if its set of isometry classes of cells is finite, and its induced piecewise Euclidean metric is a CAT(0) metric (as defined in \cite{Bridson1999}*{Def I.7.4}).
\end{definition}

Note that in general, the induced piecewise Euclidean (pseudo)metric on a polyhedral complex may fail to be complete or even geodesic. In fact, it may fail to be a metric at all. However, these pathologies are avoided when the cell complex has a finite set of isometry classes of cells \cite{Bridson1999}*{Thm~I.7.19}. We denote the set of isometry classes of cells in a Coxeter cell complex $X$ by $\op{Shapes}(X)$.

\begin{example}
    Let $W_\Gamma$ be a Coxeter group. The \textbf{Davis complex} (or Davis--Moussong complex) $X=X_\Gamma$ has a Coxeter cell of type $\Gamma'$ for each coset $gW_{\Gamma'}$, where $\Gamma'\subseteq \Gamma$ is the defining graph of a finite Coxeter subgroup~$W_{\Gamma'}\leq W_\Gamma$. The faces of Coxeter cells in $X$ are identified whenever they correspond to the same coset. 
    Note that $W_\Gamma$ acts freely and transitively on the vertices of $X$. 
\end{example}

\begin{example} If $\Gamma$ is the complete graph on $n$ vertices with edge labels all equal to 2, then $W_\Gamma=(\Z/2\Z)^n$ and $Z_\Gamma$ is an $n$-dimensional cube. Thus cube complexes are also examples of Coxeter cell complexes.
\end{example}

The work of Moussong in \cite{Moussong1988} provides a combinatorial condition on the links of vertices in~$X$ that guarantees that $X$ is CAT(0).
The \textbf{combinatorial link} of a vertex $v\in X^{(0)}$ has a $n$-simplex for each $(n+1)$-cell incident to $v$, where each face of the simplex corresponding to a cell $Z\subset X$ corresponds to a lower dimensional cell of $Z$ incident to $v$. Vertices of the combinatorial link correspond to edges of $X$ and these vertices are connected by an edge whenever the corresponding edges in $X$ bound a $2n$-gon. The \textbf{labeled link}, denoted by $\op{Lk}_\ell(v, X)$ is the combinatorial link of~$v$, where edges of $\op{Lk}_\ell(v, X)$ corresponding to $2n$-gons in $X$ have label~$n$.

If $X_\Gamma$ is the Davis complex of the Coxeter group $W_\Gamma$, then $\op{Lk}_\ell(v, X)$ has as 1-skeleton $\Gamma$, and has an $n$-simplex for each full subgraph $\Gamma'\subset \Gamma$ with $W_{\Gamma'}$ finite. In \cite{Moussong1988}*{Cor 10.2}, Moussong shows that the \textbf{metric link} of the identity in $X_\Gamma$ is CAT(1).
As a consequence of Moussong's work and \cite[Thm II.5.4]{Bridson1999}, we obtain:

\begin{theorem}
    \label{cor:link_condition}
    Let $X$ be a Coxeter cell complex with $\op{Shapes}(X)$ finite. Then $X$ is non-positively curved if and only if every vertex in $X$ is \textbf{labelled-flag}: 
    \begin{enumerate}
        \item $\op{Lk}_\ell(v, X)$ is a simplicial complex,
        \item every full subgraph $\Gamma\subset \op{Lk}_\ell(v, X)$ with $W_\Gamma$ finite spans a simplex.
    \end{enumerate}
\end{theorem}

When $X$ is a cube complex, then a full subgraph $\Gamma\subset \op{Lk}_\ell(v, X)$ has $W_\Gamma$ finite if and only if $\Gamma$ is a clique. Simplicial complexes for in which every clique span a simplex are called \textbf{flag} complexes. As such, \cref{cor:link_condition} is a natural generalization of the flag link condition typically given as the defining property of non-positively curved cube complexes \cite[Thm II.5.20]{Bridson1999}.

\subsection{Walls in CAT(0) Coxeter cell complexes}\label{sec:walls}

Here, we develop the basic theory of walls in CAT(0) Coxeter cell complexes. As we will show, many of the definitions and properties in this section bear a striking resemblance to those present in both the theories of Coxeter groups and CAT(0) cube complexes. Most of the results in this section are direct generalizations of those in \cite{Kristof-Tessier25}, and so their proofs will be omitted.  We reiterate that throughout $X$ will be assumed to be a Coxeter cell complex built out of finitely many types of Coxeter cells. In particular, $X$ is a complete geodesic metric space.

A \textbf{midcell} in a Coxeter cell $Z=Z_\Gamma$ is the fixed point set of a reflection in $W_\Gamma$. Note that each edge $e$ in $Z$ intersects a unique midcell in $Z$, denoted by $M(e,Z)$. We say two edges $e,e'\in Z$ are equivalent if $M(e,Z)=M(e',Z)$. We may then extend this equivalence relation to all edges of $X$.

\begin{definition}
        The \textbf{wall} in $X$ dual to an edge $e$ of $X$ is the abstract complex $\bigsqcup_{e'\in [e]} M(e',Z)/\!\sim$, 
        where the cell $M(e', F)$ is naturally identified with the appropriate subcomplex of $M(e', Z)$ whenever $F$ is a face of $Z$.
\end{definition}

\begin{figure}
    \centering
    \includegraphics[width=6cm]{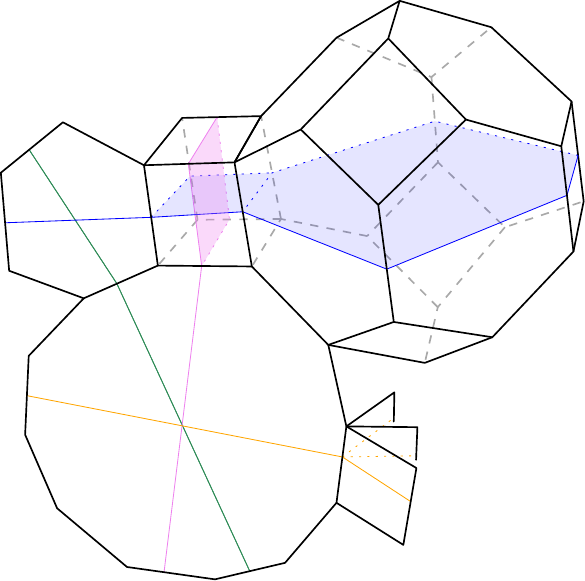}
    \caption{Examples of walls in a CAT(0) Coxeter cell complex}
    \label{fig:walls_example}
\end{figure}

\begin{remark}
    In the case where $X$ is the Davis complex of a Coxeter group, walls are fixed point sets of reflections. Therefore, they are closed, convex \cite[Prop II.6.2]{Bridson1999} and separate $X$ into two connected components. Furthermore, the combinatorial distance between vertices in $X$ (defined by the path metric on the 1-skeleton of $X$) is exactly the number of walls separating them. Note that two walls in a Davis complex intersect if and only if the underlying reflections generate a finite dihedral subgroup. 
\end{remark}

\begin{proposition}[\cite{Kristof-Tessier25}*{Prop 3.6}]
\label[proposition]{prop:walls_separate} 
    Let $X$ be a $\mathrm{CAT(0)}$ Coxeter cell complex, and $\W$ a wall in $X$. Then 
    \begin{enumerate}
        \item The map $f:\W\to X$ that sends cells in $\W$ to their appropriate midcells in $X$ is an isometric embedding. In particular, $\W$ can be identified with a closed and convex subset of $X$.
        \item Under this identification, $N_{\frac{1}{2}}(\W)$ convex, and isometric to $\W \times \left(-\frac{1}{2},\frac{1}{2}\right)$,
        \item $\W\subset X$ separates $X$ into two convex components, called (open) \textbf{halfspaces}.
    \end{enumerate}
\end{proposition}

In what follows, the wall $\W$ will always be identified with its image $f(\W)\subset X$. 
  \begin{corollary}[\cite{Kristof-Tessier25}*{Cor 3.7}]\label[corollary]{cor:isom_embedding}
    The combinatorial distance (denoted by $d_1$) between vertices $v,w\in X^{(0)}$ is equal to the number of walls in $S(v,w)$, where $S(v,w)$ is the set of walls separating $v$ from $w$.
\end{corollary}

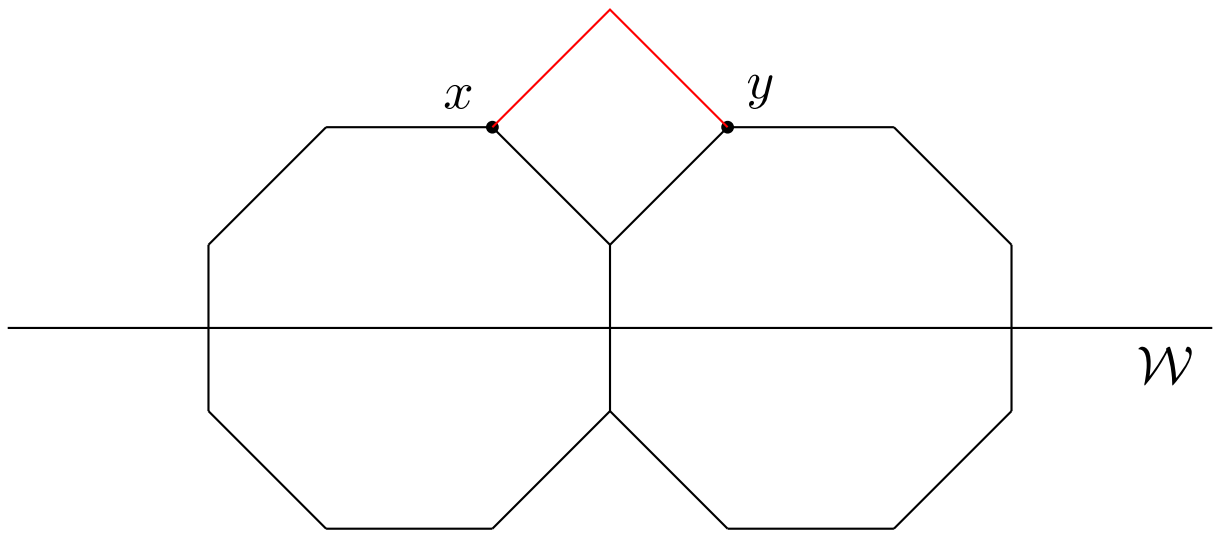
\begin{figure}
    \centering
    \centering
\begin{tikzpicture}
\def\len{1.1}
\def\nq{4}
\def\nl{((2*\nq)/(\nq-2))}

\def\rq{0.5*\len cm / sin{(180/(2*\nq))}}
\def\Rq{0.5*\len cm / tan{(180/(2*\nq))}}

\def\rl{0.5*\len cm / sin{(180/(2*\nl))}}
\def\Rl{0.5*\len cm / tan{(180/(2*\nl))}}

\foreach \k in {1,2,...,7} {
    \coordinate (v\k) at ($({180/(2*\nq)+(\k-1)*180/(\nq)}:\rq) - (\Rq, 0)$);
    \coordinate (q) at ($({180/(2*\nq)+\k*180/(\nq)}:\rq) - (\Rq, 0)$);
    
    \draw (v\k) -- (q);
}

\foreach \k in {0,1,...,7} {
    \coordinate (u\k) at ($({180/(2*\nl)+(\k-1)*180/(\nl)}:\rl) + (\Rl, 0)$);
    \coordinate (q') at ($({180/(2*\nl)+\k*180/(\nl)}:\rl) + (\Rl, 0)$);
    
    \draw (u\k) -- (q');
}

\filldraw (v2) circle (1pt) node [above left] {$x$};
\filldraw (u3) circle (1pt) node [above right] {$y$};

\draw [red] (v2) -- ($(v2)!1!90:(v1)$) -- (u3);

\draw (-3*\Rq, 0) -- (3*\Rl, 0) node [below left] {$\W$};

\end{tikzpicture}
    \unskip
    \caption{$N(\W)$ fails to be convex in both the CAT(0) metric and the combinatorial metric.}
    \label{fig:carriers_not_convex}
\end{figure}

The \textbf{carrier} of a wall $\W$, denoted by $N(\W)$ is defined to be the union of all cells in $X$ intersecting~$\W$.  In contrast to cube complexes, the carrier is in general distinct from $N_{\frac{1}{2}}(\W)$, and may fail to be convex both with the CAT(0) metric and the combinatorial metric on $X^{(1)}$ (\cref{fig:carriers_not_convex}). However, the inclusion $N(\W)^{(1)}\hookrightarrow X^{(1)}$ remains an isometric embedding \cite{Kristof-Tessier25}*{Lem 3.10}.

Many of the properties enjoyed by walls are inherited by intersections of walls. 

\begin{definition}[c.f. \cite{HruskaWise14}*{Def 2.5}]
    Let $X$ be a CAT(0) Coxeter cell complex. A \textbf{codimension-$k$ wall} $\W^k$ is an intersection of $k$ walls $\W_1,\dots,\W_{k}$ in $X$ such that for any cell $Z$ intersecting $\W^k$, $\W^k\cap Z$ is a codimension-$k$ subspace of $Z$.
\end{definition}

In particular, a codim-$1$ wall is a wall as previously defined. We will also consider $X$ to be a codim-0 wall. Note that if $\W$ is a codim-$k$ wall, then $\W$ is still closed and convex in $X$ and $N_{\frac{1}{2}}(\W)$ is isometric to $\W\times B^k$, where $B^k$ is the ball of radius $\frac{1}{2}$ in $\R^k$. 

\begin{lemma}
    \label[lemma]{lem:geom_walls}
    Let $\W^k$ be a codim-$k$ wall, and $\V$ a wall such that $\mathcal{U}=\W^k\cap \V$ is a non-empty, proper subset of $\W^k$. Then $\mathcal{U}$ is a codim-$(k+1)$ wall, and separates $\W^k$ into two connected components.
\end{lemma}

\begin{proof}
Let $\V^{\pm}$ be the two halfspaces of $\V$. Let $Z$ be a cell whose interior intersects $\mathcal{U}$. Then $(\W^k\cap \V) \cap Z \subsetneq \W^k\cap Z$, so $\mathcal{U}$ is codim-$(k+1)$ in $Z$ and $\W^k\cap \V^\pm$ are non-empty. 
Since $\V^{\pm}$ and~$\W^{k}$ are convex, the sets $\V^\pm\cap \W^{k}$ are also convex and therefore connected. It follows that $\V\cap \W^k$ separates $\W^k$ into two components.
\end{proof}

\section{Cubulating Coxeter cell complexes}\label{sec:Sageev_construction}

We briefly summarize the description of Sageev's construction by Hruska and Wise in \cite{HruskaWise14}, which in our setting associates to each CAT(0) Coxeter cell complex $X$ a connected, CAT(0) cube complex $C(X)$:

\begin{definition}[c.f. \cite{Guralnik05}*{Def 2.3, 2.6}]
    Let $\mathscr{W}$ be a (possibly infinite) set of walls in $X$. A \textbf{coorientation} of walls in $\mathscr{W}$ is an map $\sigma$ that assigns to each wall $\W\in \mathscr{W}$ one of its halfspaces~$\sigma(\W)$. A coorientation is \textbf{consistent} if $\sigma(\W_1)\cap \sigma(\W_2)\neq \emptyset$ for any pair $\W_1,\W_2\in \mathscr{W}$. We say a coorientation satisfies the \textbf{Descending Chain Condition} (DCC) if every sequence of halfspaces $\sigma(\W_1)\supseteq \sigma(\W_2)\supseteq~\dots$ with $\W_k\in \mathscr{W}$ has a minimal halfspace $\sigma(\W_n)$.
\end{definition}

The vertex set of $C(X)$ is the set of all consistent coorientations of walls in $X$ satisfying the DCC. We put an edge between two coorientations if and only if the two maps differ on a single wall of $X$. Finally, vertices in $C(X)$ span a cube if and only if they span the 1-skeleton of a cube in $C(X)^{(1)}$.
There is a natural bijection between the hyperplanes of $C(X)$ and the walls in $X$. By \cref{cor:isom_embedding}, we have an isometric embedding $X^{(1)}\subset C(X)^{(1)}$.

Observe that the edges in $C(X)$ incident to a vertex $x\in X^{(0)}\subset C(X)^{(0)}$ are in 1-1 correspondence with the set of walls that are not separated from $x$ by any other wall in $X$. We call this set the set of \textbf{$x$-elementary walls}, and denote it by $\mathcal{E}_x$. If $X$ is locally finite and satisfies the PW Property, then this set is finite for any $x\in X^{(0)}$. However, this does not a priori give us local finiteness for all of $C(X)$.
More precisely, the following is a special case of \cite[Thm 5.7]{HruskaWise14} (see also \cite{NibloReeves2003}*{Thm~3}):

\begin{theorem}\label{thm:Hruska-Wise}
    Suppose $X$ is locally finite. Then $C(X)$ is locally finite if and only if $X$ has \textbf{compact-wall separation}: there is a constant $f(K)$ for each compact $K\subset X$ such that any wall at distance $\geq f(K)$ from a compact subset $K$ is separated from $K$ by another wall.
\end{theorem}

Under the assumption that $X$ is equipped with a geometric group action, it can be deduced from \cite{Guralnik06}*{Prop 3.7, Thm B} that the PW Property implies compact-wall separation. In general, we will need the following lemma:

\begin{lemma}\label[lemma]{lem:geods_crossing_elem_walls}
    Suppose that $X$ satisfies the PW Property. Then there is a constant $C_0$ such that for any vertices $x,y\in X$ and $\W\in \mathcal{E}_x$ separating $x$ from $y$, there is a vertex $z$ lying on a combinatorial geodesic from $x$ to $y$ separated from $x$ by $\W$ at combinatorial distance $\leq C_0$ from $x$.
\end{lemma}

\begin{figure}
    \centering
    \begin{tikzpicture}
  \def\l{1.5}

  \draw (-2.5*\l,0) -- (2.5*\l,0) node [right] {$\mathcal{W}$};

  \coordinate (z) at (0,-0.5*\l);

  \draw (z) -- (0,0.5*\l);
  \filldraw (z) circle (1pt) node[below] {$z$};

  \def\rHex{0.5*\l cm / tan{(180/6)}}
  \def\RHex{0.5*\l cm/ sin{(180/6)}}

  \foreach \k in {0,1,...,5} {
    \draw ($(\k*60+30:\RHex) - (\rHex, 0)$) coordinate (lHex\k) -- ($({(\k+1)*60+30}:\RHex) - (\rHex, 0)$);
    \draw ($(\k*60+30:\RHex) + (\rHex, 0)$) coordinate (rHex\k) -- ($({(\k+1)*60+30}:\RHex) + (\rHex, 0)$);
  }

  \draw[blue, shorten <= -1.25*\l cm, shorten >= -1.5*\l cm] ($(lHex4)!.5!(lHex5)$) -- ($(lHex1)!.5!(lHex2)$);

  \draw[blue, shorten <= -1.25*\l cm, shorten >= -1.5*\l cm] ($(rHex3)!.5!(rHex4)$) -- ($(rHex1)!.5!(rHex0)$);

  \draw[thick, red, name path=midcell] ($(rHex2)!.5!(rHex3)$) -- ($(rHex0)!.5!(rHex5)$);

  \coordinate (x) at (0.8*\l, 1.5*\l);
  \filldraw (x) circle (1pt) node[above right] {$x$};
  \draw[name path=geod] (z) -- (x);
  \filldraw[red, name intersections={of=geod and midcell}] (intersection-1) coordinate (m) circle (1pt) node [below right] {$m$};

  \path (x) -- (m) -- (2.5*\l,0) coordinate (W)
    pic [draw=gray, "\tiny{\color{gray}$\geq\!\alpha$}", right, angle radius=6*\l, angle eccentricity=1] {angle=W--m--x};
\end{tikzpicture}
    \unskip
    \vspace{3.5em}
    \caption{The geodesic $[z,x]$ passes through the union of cells in $N(\W)$ containing $z$.}
    \label{fig:geods_crossing_elem_walls}
\end{figure}
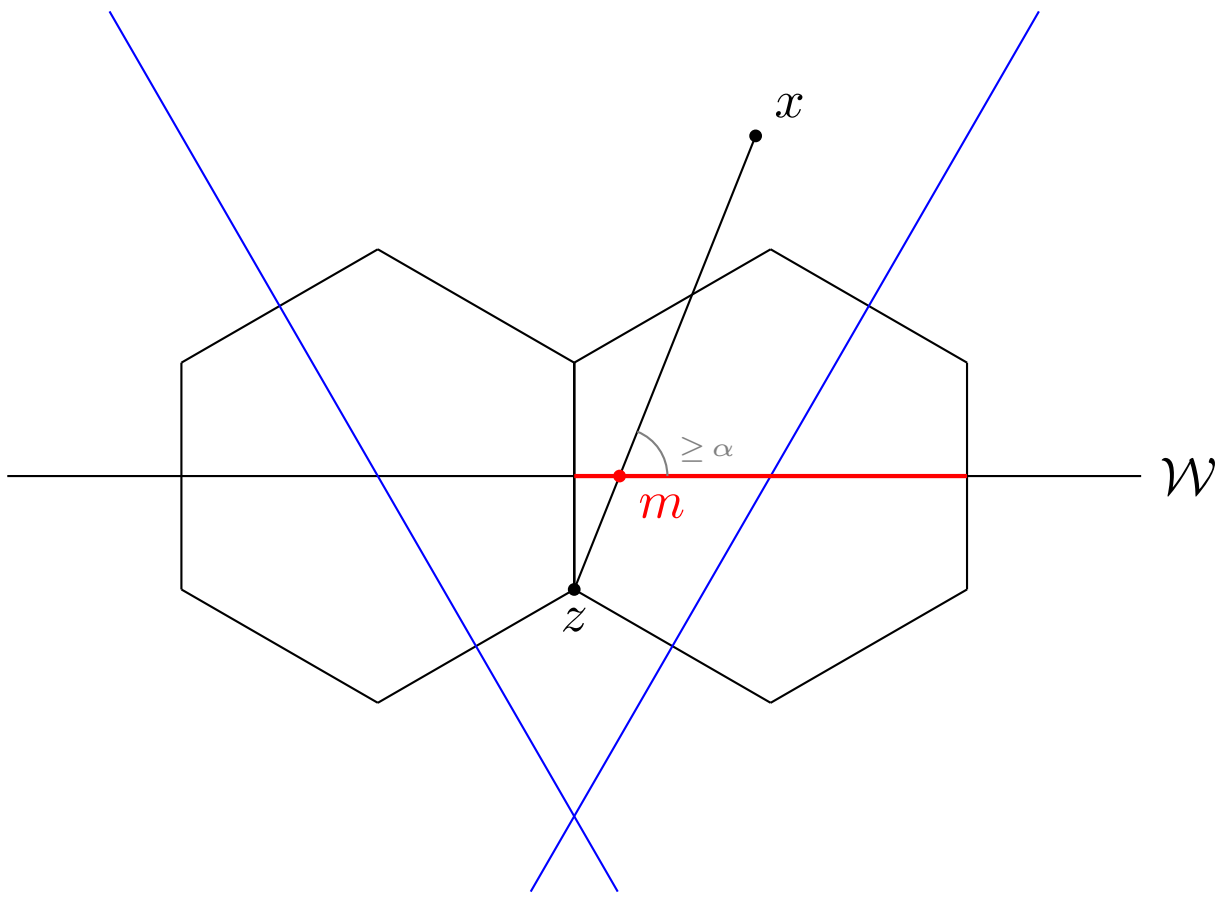

\begin{proof}
    The proof is identical to that of \cite{osajdaPrzytycki2022coxeter}*{Lem 4.3}. We reformulate their proof in the language of CAT(0) Coxeter cell complexes. Let $z\in X^{(1)}$ be a vertex as desired at minimal combinatorial distance from $x$. Note that by minimality of $z$, $\W$ is the only wall adjacent to $z$ separating $x$ from $z$. 
    
    Since $X$ is quasi-isometric to $X^{(1)}$, it suffices to uniformly bound the CAT(0) distance $d(z,x)$. We use the following quantities to do so:
    \begin{itemize}
        \item Let $\alpha$ be the minimal angle formed by a geodesic $[m, m']$ in the midcell of a cell $Z\subset X$ and a geodesic connecting $m$ to a vertex of $Z$.
        \item Let $C_1$ be the maximal distance between a vertex in a cell $Z$ and a point in a midcell of $Z$.
        \item Let $C_2$ be the uniform bound for $d(x,\pi_\W(x))$ guaranteed by the PW Property.
    \end{itemize} 

    Now, let $m$ be the intersection of the CAT(0) geodesic $[z,x]$ with $\W$. Since $\W$ is the unique wall adjacent to $z$ intersecting $[z,x]$, the segment $[z,m]$ must be contained in a cell $Z$ (see \cref{fig:geods_crossing_elem_walls}). Let $Z$ be the unique such cell containing $m$ in its interior. Then we see that $\angle_m(x,\pi_\W(x))\geq \alpha$, and so by the CAT(0) inequality,
    \[
        d(z,x)\leq C_1 + d(m,x) \leq C_1 + C_2\sin\alpha,
    \]
    giving us the desired bound.
\end{proof}

\begin{proposition}\label[proposition]{prop:cpt-wall_separation}
    Suppose $X$ satisfies the PW Property. Then there are constants $(C_n)_{n\in\N}$ such that any vertex $x\in X$ at combinatorial distance $\geq C_n$ from a wall $\W$ is separated from $\W$ by at least $n$ distinct walls. 
\end{proposition}

\begin{proof}
    We claim that we can take $C_n=n\max\{C,C_0\}$, where $C$ and $C_0$ are the constants from the PW Property and \cref{lem:geods_crossing_elem_walls} respectively. Indeed, suppose that $x$ is at distance $\geq C_n$ from a wall~$\W$. By the PW Property, there is an $x$-elementary wall $\W'$ separating $x$ from $\W$. By \cref{lem:geods_crossing_elem_walls}, there is a vertex $x'$ lying on a combinatorial geodesic from $x$ to $\W$ of length $d(x,\W)$ such that $\W'$ separates $x'$ from $x$, and $d(x', \W)\geq C_n-C_0\geq C_{n-1}$. By induction, $x'$ is separated from $\W$ by at least $n-1$ walls. But such walls must also separate $x$ from $\W$, and so $x$ is separated from $\W$ by at least $n$ walls.
\end{proof}

\begin{corollary}\label[corollary]{cor:loc_fn_cubulations_2}
    If $X$ is locally finite and satisfies the PW Property, then $C(X)$ is locally finite.
\end{corollary}

\begin{proof}
    If $X$ is locally finite, then the number of walls that can intersect a given compact set is also bounded. Therefore, \cref{prop:cpt-wall_separation}
    implies that $X$ has \textit{ball-wall separation} (i.e. compact-wall separation, where the value of $f(K)$ depends only on $\op{diam}(K)$). Thus $C(X)$ is locally finite by \cref{thm:Hruska-Wise}.
\end{proof}

\section{Angles between walls}\label{sec:angles}

In analogy to Coxeter groups, we have a well-defined notion of angle between two intersecting codim-$k$ walls $\W_1$ and $\W_2$, provided they are contained in a common codim-$(k-1)$ wall $\V$.
Fix first a cell $Z$ of $\V$ whose interior intersects $\W_1\cap \W_2$. Then $Z\cap\W_1$ and $Z\cap\W_2$ correspond to Euclidean hyperplanes in $Z$ intersecting at the center of $C$. There are now two choices of normal vectors for each hyperplane which, by \cref{lem:geom_walls}, are determined by a coorientation $\sigma$ of $\W_1$ and $\W_2$ in $\V$. The \textbf{angle between $\W_1$ and $\W_2$} in $\V$ with respect to $\sigma$ is $\pi$ minus the angle between the associated normal vectors in $Z$, and is denoted by $\angle_{\V,\sigma}(\W_1,\W_2)$. 
Since this quantity does not change when restricting to a face of $Z$ containing $\W_1\cap\W_2$, the angle $\angle_{\V,\sigma}(\W_1,\W_2)$ does not depend on our choice of $Z$.

\begin{remark}\label[remark]{rmk:finite_angles}
    If $\W_1,\W_2$ and $\W_3$ are codim-$k$-walls contained in a common codim-$(k-1)$ wall $\V$ and $\W_1\cap\W_2\cap\W_3$ is a codim-$(k+2)$ wall, then the angle between $\W_1\cap\W_2$ and $\W_1\cap\W_3$ can be written in terms of the angles between $\W_1, \W_2$ and $\W_3$. Indeed, let $\theta_{ij}=\angle_{\V,\sigma}(\W_i,\W_j)$ and let $\alpha=\angle_{\W_1,\sigma}(\W_1\cap\W_2, \W_1\cap\W_3)$. Then by the spherical law of  cosines,
    \[
        \cos(\alpha) = \frac{\cos\theta_{23}-\cos\theta_{12}\cos\theta_{13}}{\sin\theta_{12}\sin\theta_{13}}
    \]
    As such, if $X$ has finitely many types of Coxeter cells, the set of values taken by angles between codim-$k$ walls is finite.
\end{remark}

\begin{remark}
    Suppose that $X$ is the Davis complex of a Coxeter group $W$ and that $\W_r,\W_q$ are intersecting walls fixed by reflections $r,q\in W$. Let $\sigma$ be the coorientation of the $\W_r$ and $\W_q$ containing the identity in $W$. Then $\cos(\angle_{X,\sigma}(\W_r, \W_q))=-\<\alpha_r, \alpha_q\>$, where $\alpha_r, \alpha_q$ are the positive roots associated to $r$ and $q$ in the Tits representation of $W$. 
\end{remark}

\begin{proposition}\label[proposition]{prop:angle_ineq}
    Let $\W_1, \W_2$ be intersecting codim-$k$ walls lying in a codim-$(k-1)$ wall $\V$, and let $\gamma_1=[x,z]$, $\gamma_2=[z,y]$  be geodesics in $\W_1,\W_2$ respectively, such that the union $[x,z]\cup [z,y]$ is a geodesic in $\W_1\cup\W_2$, with the intrinsic path metric. Then 
    \[
            \angle(\gamma_1,\gamma_2)\geq \angle_{\V,\sigma}(\W_1,\W_2),
    \]
    where $\sigma$ is the coorientation such that $\sigma(\W_2)$ contains $\gamma_1$ and $\sigma(\W_1)$ contains $\gamma_2$.
\end{proposition}

We first prove \cref{prop:angle_ineq} for arbitrary hyperplanes in $\R^3$. Here, we reuse the notation of $\angle_{\R^3, \sigma}$ to refer to the angle between cooriented hyperplanes in $\R^3$.

\begin{lemma}
\label[lemma]{lem:angle_ineq_euc} Let $H_1, H_2\in \R^3$ be distinct intersecting hyperplanes and let $\sigma$ be a coorientation of $H_1$ and~$H_2$. Let $x\in H_1 \cap \sigma(H_2), y\in H_2\cap \sigma(H_1)$, and let $z\in H_1\cap H_2$ minimize $d(x, z)+d(z, y)$.
Then $\angle_z(x, y)\geq \angle_{\R^3,\sigma}(H_1, H_2)$.
\end{lemma}

\begin{figure}
    \centering
    \begin{tikzpicture}
    \usetikzlibrary{3d, intersections}
    \def\u{3}
    \def\x{4}
    \def\y{5}
    \def\z{2.5}

    \draw[->] (0,-\z,0) -- (0,1.2*\z,0) node[right] {$Z$};
    \coordinate (O) at (0,0,0);

    \coordinate (Hxy) at ($\x*({sin(-20)},0,{cos(-20)})$);
    \coordinate (H'xy) at ($\y*({sin(60)},0,{cos(60)})$);

    \path (Hxy) -- (O) -- (H'xy)
        pic[draw, ultra thick, white, angle radius=2*\u] {angle=Hxy--O--H'xy};
    \path (Hxy) -- (O) -- (H'xy)
        pic["$\alpha$", right, draw, angle radius=2*\u, angle eccentricity=1.2] {angle=Hxy--O--H'xy};

    \draw (Hxy) -- (O) -- (H'xy);
    \draw (0,\z,0) -- ($(Hxy) + (0,\z,0)$) node[left] {$H$}
        -- ($(Hxy) - (0,\z,0)$) --  (0,-\z,0) ;
    \draw (0,\z,0) -- ($(H'xy) + (0,\z,0)$) node[right] {$H'$}
        -- ($(H'xy) - (0,\z,0)$) --  (0,-\z,0);
        
    \def\h{1.5}
    \coordinate (px) at ($0.6*(Hxy)$);
    \coordinate (py) at ($0.6*(H'xy)$);
    \coordinate (x) at ($(px)-(0,\h,0)$);
    \coordinate (y) at ($(py)+(0,\h,0)$);

    \path[name path=xy] (x) -- (y);
    \path[name path=pxpy] (px) -- (py);
    \filldraw[name intersections={of=xy and pxpy}] (intersection-1) node(z) {} circle (0pt);

    \path (y) -- (py) -- (z) pic[draw, ultra thick, white, angle radius=2*\u]
        {right angle=y--py--z};

    \draw[dashed] (O) -- (px) node[pos=0.5, above, gray] {$\lambda$} -- (x) node[pos=0.5, left] {$h$} -- cycle;
    \draw (O) -- (y) -- (py) node[pos=0.5, right] {$h$};

    \draw[ultra thick, white] (px) -- (py);
    \draw[ultra thick, white] (x) -- (y);
    \draw (px) -- (py);
    \draw[shorten >=-5pt] (y) -- (z);
    \draw[dashed, shorten >=-5pt] (x) -- (z);

    \path (Hxy) -- (O) -- (H'xy)
        pic[draw, gray, angle radius=2*\u] {right angle=y--py--z};

    \path (Hxy) -- (O) -- (H'xy)
        pic[draw, gray, angle radius=2*\u] {right angle=x--px--z};
    
    \fill[opacity=.15] (O) -- (H'xy) -- ($(H'xy)+(Hxy)$) node[opacity=1, right=0.8] {$XY$} -- (Hxy) -- cycle;

    \filldraw (O) circle (1pt) node[above left] {$\vec{0}$};
    \filldraw (x) circle (1pt) node[below left] {$x$};
    \filldraw (y) circle (1pt) node[right] {$y$};
    \filldraw (px) circle (1pt) node[above] {$\pi(x)$};
    \filldraw (py) circle (1pt) node[below] {$\pi(y)$};
    \filldraw (z) circle (.5pt);

\end{tikzpicture}
    \unskip
    \caption{Proof of \cref{lem:angle_ineq_euc}}
    \label{fig:angle_ineq_euc}
\end{figure}

\begin{proof}
    Assume without loss of generality that $z=(0,0,0)$, that $\|x\|=\|y\|$, and that $H\cap H'$ is the $Z$-axis. Let $\pi(x)$ and $\pi(y)$ be the projections of $x$ and $y$ onto the $XY$-plane. Note that $\angle_z(\pi(x), \pi(y))=\angle_{\R^3,\sigma}(H_1, H_2)=: \alpha$. 
    Write $x=\pi(x)+h(0,0,1)$ and $y=\pi(y)+h'(0,0,1)$. By our assumption, the function
    \[
            f(t) = \|x-t(0,0,1)\|+\|y-t(0,0,1)\|
    \]
    is minimized at $t=0$. Thus we have $0 = f'(0) = \frac{h}{\|x\|} + \frac{h'}{\|y\|}$, and so $h'=-h$. Now, letting $\lambda = \|\pi(x)\|= \|\pi(y)\|$ and $\beta=\angle_z(x,y)$, we see that
   \begin{equation*}
        \|\pi(x)-\pi(y)\|^2=2\lambda^2 - 2\lambda^2\cos\alpha
        \qquad\text{and}\qquad\|x-y\|^2=2(\lambda^2+h^2) - 2(\lambda^2 + h^2)\cos\beta
   \end{equation*}
    Additionally, the Pythagorean theorem gives us the relation
        $\|x-y\|^2 = \|\pi(x)-\pi(y)\|^2 + 4h^2.$
    Combining the three equalities above gives us that 
    \[
            cos\beta  =\frac{\lambda^2\cos\alpha - h^2}{\lambda^2 + h^2}
    \]
    and so $\cos\beta-\cos\alpha = \frac{h^2(1-\cos\alpha)}{\lambda^2+h^2}\geq 0$.
\end{proof}

\begin{proof}[Proof of \cref{prop:angle_ineq}]
    Set $Y=\W_1\cap \W_2$, so that $z\in Y$, and assume without loss of generality that $x,y$ are contained in $\W_1\cap N_{\frac{1}{2}}(\W_2)$ and $\W_2\cap N_{\frac{1}{2}}(\W_1)$ respectively.
Note that we have an isometric embedding $i:N_{\frac{1}{2}}(Y)\hookrightarrow Y\times \R^2$, under which $\W_i\cap N_{\frac{1}{2}}(\W_j)$ gets mapped into $Y\times L_i$ for some distinct intersecting lines $L_1$ and $L_2$ in $\R^2$. 
Since the projection of a geodesic onto a factor is itself geodesic,
$J:=\pi_Y([x,z]\cup [z,y])$ is isometric to either a line segment or a point. This gives us a further isometric embedding of $J\times \R^2\hookrightarrow \R\times \R^2$, under which $i([x, z])$ and $i([z,y])$ map into $H_1=\R\times L_1$ and $H_2=\R\times L_2$ respectively. Now, $H_1$ and $H_2$ are hyperplanes in $\R^3$, with coorientations determined by $\sigma$, and $\angle_{\R^3, \sigma}(H_1, H_2)=\angle_{\V,\sigma}(\W_1, \W_2)$. Therefore, by \cref{lem:angle_ineq_euc} we have that $\angle_z(x,y)\geq \angle_{\R^3, \sigma}(H_1, H_2)= \angle_{\V,\sigma}(\W_1, \W_2)$.
\end{proof}

\begin{corollary}
    [CAT(0) inequality for walls] \label[corollary]{cor:wall_angle_CAT(0)_ineq}
    Suppose that $\W_1,\W_2,$ and $\W_3$ are three pairwise intersecting codim-$k$ walls contained in a codim-$(k-1)$ wall $\V$. Suppose further that the walls have empty global intersection, and let $\sigma$ be the unique coorientation of walls such that $\sigma(\W_i)\supset \W_j\cap\W_\ell$ for $i,j,k$ distinct. Then
    \[
        \angle_{\V,\sigma}(\W_1, \W_2) + \angle_{\V,\sigma}(\W_1, \W_3) + \angle_{\V,\sigma}(\W_2, \W_3)  \leq \pi
    \]
\end{corollary}

\begin{proof}
    Let $Y\subset \V$ be the union $\W_1\cup\W_2\cup\W_3$, endowed with the intrinsic path metric. By \cite{Bridson1999}*{Thm II.11.1}, $Y$ is non-positively curved. Also, $Y$ has fundamental group $\Z$, generated by geodesic triangles with vertices in $\W_i\cap\W_j$. By \cite{Bridson1999_semisimple}, the action of $\Z$ on the universal cover $\widetilde{Y}$ is hyperbolic, and so there is a geodesic line in $\widetilde{Y}$ fixed by $\Z$. Projecting this line onto $Y$ yields a geodesic triangle $\Delta$ with edges $\gamma_i\subset \W_i$. Since $\tilde{Y}\to Y$ is a local isometry, $\gamma_i\cup \gamma_j$ is geodesic in $\W_i\cup \W_j$ for each $1 \leq i,j\leq 3$. By \cref{prop:angle_ineq}, we thus have that $\angle_{\V,\sigma}(\W_i, \W_j) \leq \angle(\gamma_i, \gamma_j)$. 
    Since $\V$ is CAT(0), \cref{prop:cat0_angle_ineq} yields
    \[
        \sum_{i<j}\angle_{\V,\sigma}(\W_i, \W_j) \leq \sum_{i< j}\angle(\gamma_i, \gamma_j) \leq \pi
    \]
\end{proof}

The angles between walls provide a method to easily calculate the combinatorial distance between a vertex and a given wall, generalizing \cite{BrinkHowlett1993}*{Lem 1.7}:

\begin{proposition}
    \label[proposition]{lem:angle_lemma} 
Let $u,v\in X^{(0)}$ be adjacent vertices separated by a wall $\W'$, and let $\W$ be a wall intersecting $\W'$. Take $\sigma$ to be the coorientation of $\W$ and $\W'$ containing $u$. Then
\[
    d_1(u, \W) = \begin{cases}
        d_1(v, \W) + 1 & \text{if } \angle_{X,\sigma}(\W, \W') > \pi/2 \\
        d_1(v, \W) & \text{if } \angle_{X,\sigma}(\W, \W')= \pi/2 \\
        d_1(v, \W) - 1 & \text{if } \angle_{X,\sigma}(\W, \W') < \pi/2 \\
    \end{cases}
\]
\end{proposition}

The idea of the proof of \cref{lem:angle_lemma} is to show that a vertex of $\W$ close to the CAT(0) projection $\pi_\W(u)$ minimizes $d_1(u, \W)$. More precisely, we choose a particular vertex $p$ lying in the same closed halfspaces as~$\pi_\W(m)$, where $m$ is the midpoint of the edge $[u,v]$. Note that there are many choices for such~$p$, specified by the choice of halfspaces of the walls containing $\pi_\W(m)$. We then use \cref{prop:angle_ineq} to show that the halfspace of $\W'$ containing $\pi_\W(m)$ is determined by $\angle_{X,\sigma}(\W, \W')$ (\cref{cor:angle_between_m_and_pi}).
Finally, to compare $d_1(u, p)$ with $d_1(u,\W)$, we introduce the notion of a reflection and show that each wall separating $u$ from $p$ has a reflection that does not separate $u$ from $p$, except possibly for $\W'$ in the case where $\W'$ is perpendicular to $\W$. 

\begin{definition}
    Suppose that $\V$ and $\W$ are intersecting walls in a CAT(0) Coxeter cell complex $X$. Let $Z_\Gamma$ be a cell whose interior intersects $\V\cap \W$. Then the intersection of $\V$ and $\W$ with $Z_\Gamma$ are midcells in $Z_\Gamma$. Thus $\W\cap Z_\Gamma$ is the fixed point set of a (unique) reflection $r\in W_\Gamma$. The midcell $r\cdot(\V\cap Z_\Gamma)$ is the cell of a new wall $\overline{\V}$ in $X$ called the \textbf{reflection of $\V$ along $\W$}. One may check that the wall obtained through this procedure does not depend on our choice of $Z$. 
\end{definition}

\begin{lemma}
    \label[lemma]{cor:angle_between_m_and_pi}
    Let $x\in X$, let $\W$ be a wall that does not contain $x$, and
    let $p$ be a vertex of $\W$ lying in all the same closed halfspaces as $\pi_\W(x)$. Suppose $\V$ is a wall intersecting $\W$ that either 
    \begin{enumerate}
        \item separates $x$ from $p$, or
        \item contains $x$ but not $\pi_\W(x)$.
    \end{enumerate}
Then $\angle_{X,\sigma}(\V,\W)<\pi/2$, where $\sigma(\W)$ is chosen to contain $x$ and $\sigma(\V)$ is chosen to contain $p$.  Furthermore, the reflection $\overline{\V}$ of $\V$ along $\W$ does not satisfy either of the conditions (1) and (2) above.
\end{lemma}
\begin{figure}
    \centering
    \[\vcenter{\hbox{
\begin{tikzpicture}
    \usetikzlibrary{decorations.pathmorphing}
    
    \def\l{5}

    \coordinate (px) at (0.6*\l, 0);
    \coordinate (x) at ($(px)+(0, .6*\l)$);
    \coordinate (p) at (0.75*\l, 0);

    \coordinate (V) at (30:\l);
    \coordinate (W) at (\l, 0);
    
    \draw[shorten <=-5*\l ] (0,0) -- (W) node[right] {$\mathcal{W}$}; 
    \draw[shorten <=-5*\l , name path = V] (0,0) -- (V) node[above right] {$\mathcal{V}$}; 

    \filldraw (px) circle (1pt) node [below] {$\pi_\mathcal{W}(x)$};
    \filldraw (x) circle (1pt) node [above] {$x$};
    \filldraw (p) circle (1pt) node [below] {$p$};
    
    \draw[name path = geod, decorate, decoration={snake, amplitude=0.25mm, segment length=5mm}] (px) -- (x);

    \filldraw[name intersections={of=geod and V}] (intersection-1) circle (1pt) node [above left] {$x'$};

    \path (V) -- (0,0) coordinate (O) -- (W)
        pic[draw=gray, "\small{\color{gray} $<\!\pi/2$}", right, angle radius=4*\l, angle eccentricity=1.1] {angle=W--O--V};

\draw[gray, ->] (30:.95*\l) -- ($(30:.95*\l)!2!-90:(V)$);
\draw[gray, ->] (.95*\l, 0) -- ($(.95*\l, 0)!2!90:(W)$);
\end{tikzpicture}}}
\quad\quad
\vcenter{\hbox{
\begin{tikzpicture}
    \def\w{2}
    \coordinate (px) at (0,0,0);
    \coordinate (x) at (0,3,0);
    \coordinate (y) at (2,0,0);

    \draw (-\w,0,\w) -- (-\w,0,-\w) -- (1.7*\w,0,-\w) node[above] {$\mathcal{W}$} -- (1.7*\w,0,\w) -- cycle;

    \coordinate (W_arrow_base) at ($0.8*(2*\w,0,-\w)$);
    \draw[white, ultra thick] (W_arrow_base) -- ($(W_arrow_base)+(0,0.6,0)$);
    \draw[gray, ->] (W_arrow_base) -- ($(W_arrow_base)+(0,0.6,0)$);

    \coordinate (V_vect) at (xyz cylindrical cs:radius=1,angle=30);
    \coordinate (V_perp) at (xyz cylindrical cs:radius=1,angle=-60);
    \coordinate (V_arrow_base) at ($0.9*3*(V_vect) + 0.9*(0,0,-\w)$);
    \draw[->, gray, shorten <=4pt] (V_arrow_base) -- ($(V_arrow_base) + 0.6*(V_perp)$);

    \filldraw[white] ($-1*(V_vect)+(0,0,\w)$) -- ($3*(V_vect)+(0,0,\w)$) -- ($3*(V_vect) + (0,0,-\w)$) -- (0,0,-\w);
    \draw ($-1*(V_vect)+(0,0,\w)$) -- ($3*(V_vect)+(0,0,\w)$) -- ($3*(V_vect) + (0,0,-\w)$) node(V_label) [right] {$\mathcal{V}$} -- (0,0,-\w);
    \draw[shorten >= 4.5em] ($-1*(V_vect)+(0,0,\w)$) -- ($-1*(V_vect)-(0,0,\w)$);

    \filldraw (px) circle (1pt) node[left] {$\pi_\mathcal{W}(x)$};
    \filldraw (x) circle (1pt) node[left] {$x$};
    \filldraw (y) circle (1pt) node[right] {$y$};

    \draw[shorten >= 2.7em] (x) -- (2,0,0) node(q) [pos=0.73] {};
    \draw[shorten <= 8.5em] (x) -- (2,0,0);
    \filldraw (q) circle (1pt) node[above right] {$x'$};
    
    \draw[dashed] (q) -- (px) -- (y);
    \draw (x) -- (px);
    \draw (0,0, \w)-- (0,0,-\w) node {};

    \filldraw (1,0,1) circle (1pt) node[right] {$p$};

    \draw (1.7*\w,0,\w) -- (-\w,0,\w);
\end{tikzpicture}}}
\]
    \unskip
    \caption{Proof of \cref{cor:angle_between_m_and_pi}}
\end{figure}

\begin{proof}
    We distinguish two cases. If $\V$ does not contain $\pi_\W(x)$, then the segment $[x,\pi_\W(x)]$ intersects~$\V$ at some point $x'$ distinct from $\pi_\W(x)$. By the CAT(0) inequality for angles (\cref{prop:cat0_angle_ineq}), $\angle_z(x',\pi_\W(x))<\pi/2$ for any $z\in \V\cap\W$. 
    Picking $z\in \V\cap\W$ lying on the geodesic between~$x'$ and~$\pi_\W(x)$ in $\V\cup \W$ (with the intrinsic path metric on $\V\cup \W$) and applying \cref{prop:angle_ineq} then gives us that $\angle_{X,\sigma}(\V,\W)\leq \angle_z(x',\pi_\W(x))<\pi/2$. 
    Additionally, $p$ and $\pi_\W(x)$ lie in the same halfspace of $\V\cap \W = \overline\V\cap\W$ (as a codim-2 wall in $\W$), and therefore they lie in the same halfspace of~$\overline\V$. 
    So, if $\overline{\V}$ satisfied (1) or (2), then after choosing $\sigma(\overline{\V})$ to be the halfspace containing $p$, we would have that $\angle_{X,\sigma}(\overline{\V},\W)=\pi-\angle_{X,\sigma}(\V,\W) < \pi/2$, which is a contradiction.

    Now, if  $\V$ contains $\pi_\W(x)$, then $\pi_\W(x)$ is contained in the interior of some cell $Z$, and $p$ is the midpoint of an edge of $Z$. 
    We may assume without loss of generality that $x$ lies in $Z$.
    Choose $y\in \sigma(\V)\cap\W$ so that the projection of $y$ onto $\V\cap\W$ is $\pi_\W(x)$. Then $\V$ separates $x$ from $y$, and so the segment $[x,y]$ intersects $\V$ at some point $x'$. But then 
    \[
        \angle_{X,\sigma}(\V,\W)=\angle_{\pi_\W(x)}(x',y)<\angle_{\pi_\W(x)}(x,y)=\pi/2
    \]
    as desired. A similar argument as in the first case then shows that $\overline\V$ does not separate $x$ from $y$, and therefore does not separate $x$ from $p$. 
\end{proof}

\begin{figure}
    \def\l{2}
    \centering
\begin{subfigure}{0.30\textwidth}
    \centering
    \begin{tikzpicture}
        \def\a{143}
        
        \draw (-0.5*\l,0) -- (\l,0) node[right] {$\W$};
        \draw (180-\a:1.2*\l) node [above right] {$\W'$} -- (-\a:0.5*\l);
        \coordinate (tmp) at ($(0,0)!.8!(180-\a:\l)$);
        \coordinate (1) at ($(tmp)!0.5!-90:(180-\a:\l)$);
        \coordinate (s) at ($(tmp)!0.5!90:(180-\a:\l)$);
        \filldraw (s) circle (1pt)  node[above] {$u$};
        \filldraw (1) circle (1pt)  node[right] {$v$};

        \draw (180-\a:0.15*\l) arc[radius=0.15*\l, start angle=180-\a, end angle= 180];

        \draw[dotted] (1) -- (s);
        
    \end{tikzpicture}
    \end{subfigure}
\begin{subfigure}{0.30\textwidth}
    \centering
    \begin{tikzpicture}
        \def\a{90}
        
        \draw (-0.75*\l,0) -- (0.75*\l,0) node[right] {$\W$};
        \draw (180-\a:\l) node [above] {$\W'$} -- (-\a:0.3*\l);     
        \coordinate (tmp) at ($(0,0)!.8!(180-\a:\l)$);
        \coordinate (1) at ($(tmp)!0.5!-90:(180-\a:\l)$);
        \coordinate (s) at ($(tmp)!0.5!90:(180-\a:\l)$);
        \filldraw (s) circle (1pt)  node[above left] {$u$};
        \filldraw (1) circle (1pt)  node[above right] {$v$};

        \draw[dotted] (1) -- (s);

        \draw (0,0.1*\l) -- (-0.1*\l,0.1*\l) -- (-0.1*\l,0);
        
    \end{tikzpicture}
\end{subfigure}
\begin{subfigure}{0.30\textwidth}
    \centering
    \begin{tikzpicture}
        
        \def\a{37}
        
        \draw (-\l,0) -- (0.5*\l,0) node[right] {$\W$};
        \draw (180-\a:1.2*\l) node [above left] {$\W'$} -- (-\a:0.5*\l);     

        \coordinate (tmp) at ($(0,0)!.8!(180-\a:\l)$);
        \coordinate (1) at ($(tmp)!0.5!-90:(180-\a:\l)$);
        \coordinate (s) at ($(tmp)!0.5!90:(180-\a:\l)$);
        \filldraw (s) circle (1pt)  node[left] {$u$};
        \filldraw (1) circle (1pt)  node[above left] {$v$};

        \draw[dotted] (1) -- (s);

        \draw (180-\a:0.2*\l) arc[radius=0.2*\l, start angle=180-\a, end angle= 180];
        
    \end{tikzpicture}
\end{subfigure}
    \unskip
    \caption{The three cases in \cref{lem:angle_lemma}.}
    \label{fig:angle_lemma}
\end{figure}
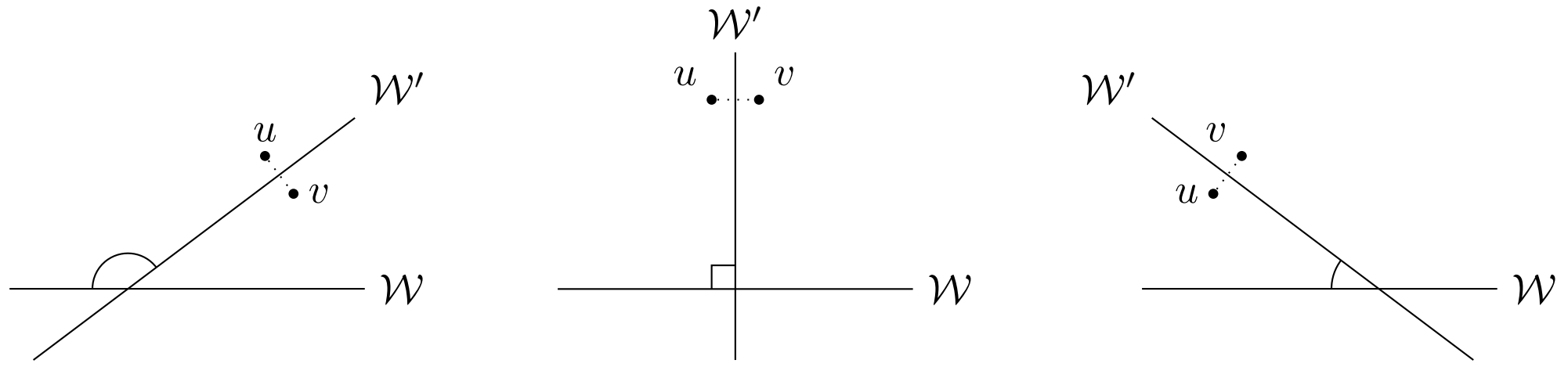

\begin{proof}[Proof of \cref{lem:angle_lemma}]
    Suppose that $p'\in \W\cap X^{(1)}$ minimizes $d_1(u, \W)$. Let $m\in \W'$ be the midpoint of the edge $[u,v]$, and let $p\in\W\cap X^{(1)}$ be the vertex of $\W$ lying in the same closed halfspaces as $\pi_\W(m)$ and separated from $p'$ by the walls containing $\pi_\W(m)$.
    
    We claim that if $\angle_{X, \sigma}(\W, \W')\neq \pi/2$, then $d_1(u,p)=d_1(u, \W)$, and $d_1(u, p)=d_1(u, \W) + 1$ otherwise.
    Observe that to compare $d_1(u,p)$ with $d_1(u,\W)$, it suffices to look at the walls in the symmetric difference $S(p,p')=S(u,p)\op{\triangle} S(u,p')$.
    Furthermore, the only wall in $S(p,p')$ that may contain both $m$ and $\pi_\W(m)$ is $\W'$ in the special case where $\angle_{X,\sigma}(\W,\W')=\pi/2$.

    \textbf{Case 1:} $\angle_{X,\sigma}(\W,\W')>\pi/2$.
    If $\V\in S(p,p')$ separates $u$ from $p$, then it either separates $m$ from $p$ or equals $\W'$. By \cref{cor:angle_between_m_and_pi}, the reflection $\overline\V$ of $\V$ along $\W$ neither separates $m$ from~$p$, nor does it contain $m$ and not $\pi_\W(m)$. Consequently, $\overline\V$ does not separate $u$ from $p$ and since $\overline\V\in S(u,p)\op{\triangle} S(u,p')$, $\overline\V$ separates $u$ from $p'$. It follows that at least half of the walls in $S(p,p')$ lie in $S(u,p)$, and thus $d_1(u,p)\leq d_1(u,p')=d_1(u,\W)$, giving us the desired claim.

    Now, if $\sigma(\W')$ contained $p$, then by \cref{cor:angle_between_m_and_pi}, we would have $\angle_{X,\sigma}(\W,\W')<\pi/2$, a contradiction. Therefore, $\W'$ separates $p$ from $u$, and so
    \[
        d_1(u,\W) = d_1(u,p) = d_1(v,p) + 1 = d_1(v,\W) + 1.
    \]

    \textbf{Case 2:} $\angle_{X,\sigma}(\W,\W')=\pi/2$. Similarly to the previous case, any wall in $S(p,p')\setminus \{\W'\}$ separating~$u$ from $p$ has a reflection that separates~$u$ from $p'$. Therefore, $d_1(u,p)= d_1(u,p')+ 1$ if~$\W'$ separates $u$ from $p$ and $d_1(u,p)= d_1(u,p')- 1$ if $\W'$ separates $u$ from  $p'$. But then by minimality of~$p'$, we must have $d_1(u, p)=d_1(u,p')+1 = d_1(u, \W)+1$. Therefore, $\W'$ separates $u$ from $p'$ and
    \[
        d_1(u,\W) = d_1(u,p') = d_1(u,p)-1 = d_1(v,p) \geq d_1(v,\W).
    \]
    The other inequality follows from interchanging $u$ and $v$, so $d_1(u,\W)=d_1(v,\W)$.

    Finally, the remaining case where $\angle_{X,\sigma}(\W,\W')<\pi/2$ follows from the first case by interchanging $u$ and $v$ as well.
\end{proof}

\section{Bipodality and the Parallel Wall theorem}\label{sec:proof}

Let $X$ be a CAT(0) Coxeter cell complex, and let $\W, \W'$ be distinct walls in $X$ sharing a common adjacent vertex $v$. The intersection of the halfspaces of $\W, \W'$ containing $v$ is called a \textbf{sector} of~$X$ bounded by $\W$ and $\W'$. A sector is said to be \textbf{sharp-angled} if its bounding walls are not perpendicular to one another.

\begin{definition}
\label{def:fund_sectors_2}
    Let $x\in X^{(0)}$ and let $\W $ and $\W'$ be two intersecting walls. Consider $\mathscr{V}$, the set of all walls containing the codim-2 wall $\W\cap \W'$, and choose $\sigma(\V)$ to be the halfspace of $\V$ containing $x$ for every $\V\in \mathscr{V}$. 
Then 
    \[
            D_x(\W, \W')=\bigcap_{\V\in \mathscr{V}} \sigma(\V)
    \]
    is a sector containing $x$, called the \textbf{canonical sector} for $\W, \W'$ containing $x$. 
\end{definition}

\begin{remark}\label[remark]{lem:fund_sectors}
    Note that $\W,\W'$ contain $\V\cap \V'$ for all $\V, \V'\in \mathscr{V}$, and that the size of $\mathscr{V}$ is bounded by the largest number of sides of a polygon in $\op{Shapes}(X)$. As such, the number of walls that may pair to give the same canonical sector is bounded by a constant depending only on $\op{Shapes}(X)$.
\end{remark}

\begin{definition}[{\cite{dyer2016smallrootslowelements}*{Def 4.14}}]
    Fix a basepoint $x$. A set of walls $\mathscr{W}$ is said to be \textbf{bipodal} if it satisfies the following: if $\W$ and $\W'$ bound a sector containing $x$, and $\V\neq \W, \W'$ is a wall containing $\W\cap \W'$, then $\V\in \mathscr{W}$ implies that $\W, \W'\in\mathscr{W}$.
\end{definition}

Recall that a wall is $x$-elementary if it is is not separated from $x$ by another wall of $X$. The definition for bipodality of a CAT(0) Coxeter cell complex $X$ given in the introduction is equivalent to requiring that the set $\mathcal{E}_x$ of $x$-elementary walls is bipodal for any basepoint $x\in X^{(0)}$. 
Note that bipodality does not hold for arbitrary CAT(0) Coxeter cell complexes (see \cref{fig:bip_counterexample}). The following is a slightly stronger version of \cref{thm:bip_BPI_implies_PW}:

\begin{figure}
    \centering
    \begin{tikzpicture}
  \def\l{0.9}
  \def\n{3}
  \def\m{6}
  \def\rHex{0.5*\l cm / tan{(180/(2*\n))}}
  \def\RHex{0.5*\l cm/ sin{(180/(2*\n))}}

  \def\rDod{0.5*\l cm / tan{(180/(2*\m))}}
  \def\RDod{0.5*\l cm/ sin{(180/(2*\m))}}

  \foreach \k in {0,1,...,5} {
    \coordinate (Hex\k) at ($(\k*60+30:\RHex)$);
    \draw[thick] (Hex\k) -- ($({(\k+1)*60+30}:\RHex)$);
  }

  \foreach \k in {0,1,...,10} {
    \coordinate (DodR\k) at ($({180/(2*\m)-(\k+7)*180/(\m)}:\RDod) + (\rHex+\rDod, 0)$);
    \coordinate (tempR) at ($({180/(2*\m)-(\k+8)*180/(\m)}:\RDod) + (\rHex+\rDod, 0)$);

    \coordinate (DodL\k) at ($({180/(2*\m)+\k*180/(\m)}:\RDod) - (\rHex+\rDod, 0)$);
    \coordinate (tempL) at ($({180/(2*\m)+(\k+1)*180/(\m)}:\RDod) - (\rHex+\rDod, 0)$);

    \draw[thick] (DodR\k) -- (tempR);
    \draw[thick] (DodL\k) -- (tempL);
  }

  \coordinate (sqR_vert) at ($(Hex1)!1!90:(Hex0)$);
  \coordinate (sqL_vert) at ($(Hex1)!1!-90:(Hex2)$);
  \draw[thick] (Hex1) -- (sqR_vert) -- (DodR1);
  \draw[thick] (Hex1) -- (sqL_vert) -- (DodL1);

  \coordinate (x) at ($(sqR_vert)!1!-120:(Hex1)$);
  \draw[thick] (sqR_vert) -- (x) -- (sqL_vert);

  \filldraw (x) circle (1pt) node[right] {$x$};
  
  \draw[red] ($(DodL5)!.5!(DodL6)$) -- ($(DodR5)!.5!(DodR6)$) node[below left] {$\mathcal{V}$};
  
  \draw[blue] ($(DodL1)!.5!(sqL_vert)$) -- ($(Hex4)!.5!(Hex5)$) node[below right] {$\mathcal{W}'$};
  \draw[blue] ($(DodR1)!.5!(sqR_vert)$) -- ($(Hex4)!.5!(Hex3)$) node[below left] {$\mathcal{W}$};
\end{tikzpicture}
    \unskip
    \caption{$\V$ is $x$-elementary, but $\W$ and $\W'$ are not.}
    \label{fig:bip_counterexample}
\end{figure}
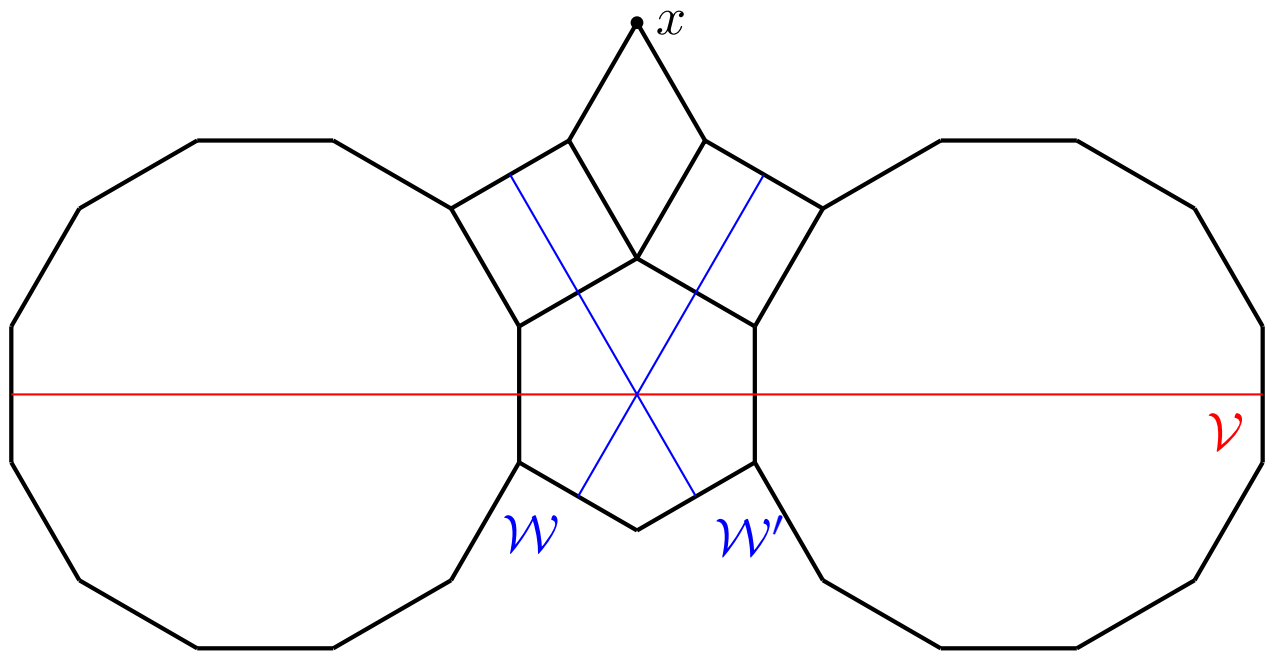

\begin{theorem}\label{thm:bipodal_par_wall_thm}
    Let $X$ be a $\mathrm{CAT(0)}$ Coxeter cell complex satisfying the Bounded Pairwise Intersection Property. 
    Then there is a bound $C$, depending only on $X$, such that for any $x\in X^{(0)}$ with $\mathcal{E}_x$ bipodal, the walls in $\mathcal{E}_x$ are at distance $\leq C$ from $x$.
\end{theorem}

\begin{proof}
    Let $x\in X^{(0)}$ be such that $\mathcal{E}_x$ is bipodal, and let $\W$ be a wall. Fix a point $p\in \W\cap X^{(1)}$ minimizing $d_1(x, \W)$, and assume that each wall in $S(x, p)$ intersects $\W$. Then for each wall in $\V\in S(x,p)$, we may consider the canonical sector containing $x$ for the pair of walls $\V$ and $\W$. Since~$\V$ is not perpendicular to $\W$ by \cref{lem:angle_lemma}, the sectors we obtain are sharp-angled, and their bounding walls are distinct from $\W$. By \cref{lem:fund_sectors}, 
    the number of distinct sectors obtained this way increases as the size of $S(x,p)$ increases. Enumerate these sectors as $D_1,\dots,D_r$. Note that one of the walls bounding $D_i$, say $\V_i$, separates $x$ from $p$.
    
    Let $N$ be the bound on the number of pairwise intersecting walls in $X$, and choose $C\geq 0$ so that $r>N$ whenever $\#S(x,p) > C$.
    Then there is a pair of disjoint walls $\V_i$ and $\V_j$ such that~$\V_i$ separates $x$ from $\V_j$. 
    But then $\V_j\notin \mathcal{E}_x$, and since $\W$  contains the intersection of the two walls bounding $D_j$, we obtain under the assumption of bipodality that $\W\notin \mathcal{E}_x$. We therefore conclude that $\W\notin \mathcal{E}_x$ for any wall $\W$ with $d_1(x,\W)>C$.
\end{proof}

\bibliography{refs}
\end{document}